\documentclass[reqno,11pt]{article}

\usepackage[a4paper, hdivide={3.0cm,,3.0cm},vdivide={3.48cm,,3.48cm}]{geometry}
\usepackage{amssymb}
\usepackage{amsmath}
\usepackage{amsthm}
\usepackage{enumitem}
\usepackage{amstext}
\usepackage{graphicx}
\usepackage{hhline}
\usepackage{psfrag}
\usepackage{float}
\usepackage[dvips]{epsfig}
\usepackage{rotating}
\usepackage{varioref}
\usepackage{color}

\usepackage[utf8]{inputenc}
\usepackage[T1]{fontenc}
\usepackage[english]{babel}
\usepackage{ngerman}
\usepackage{fancyhdr}

\usepackage{comment}

\usepackage{hyperref}

\allowdisplaybreaks

\DeclareSymbolFont{msbm}{U}{msb}{m}{n}
\DeclareMathSymbol{\N}{\mathalpha}{msbm}{'116}
\DeclareMathSymbol{\R}{\mathalpha}{msbm}{'122}

\newcommand{\interior}[1]{\overset{\circ}{#1}}

\def\Om{\Omega}

\def\pdi{\partial_i}
\def\pdj{\partial_j}

\def\emb{\hookrightarrow}

\def\sup{\text{sup}}
\def\supp{\text{supp}}

\def\N{\mathbb N }

\def\limn{\underset{n \to \infty}{\text{lim}}}

\def\lifetime{\mathcal X}
\newcommand{\mhat}[1]{\hat{#1}}

\theoremstyle{theorem} % Text italic, Heading 

\newtheorem{theorem}{Theorem}[section] % Beginn der section mit theorem, Bezeichnung Theorem, neuz{\"a}hlen nach section
\newtheorem{proposition}[theorem]{Proposition}
\newtheorem{lemma}[theorem]{Lemma}
\newtheorem{corollary}[theorem]{Corollary}

\theoremstyle{definition} % Text normal, heading fett
\newtheorem{definition}[theorem]{Definition}

\newtheorem{remark}[theorem]{Remark}
\newtheorem{condition}[theorem]{Condition}

\newcommand{\E}{\mathcal{E}}
\newcommand{\D}{\mathcal{D}}
\newcommand{\Dom}{D(\mathcal{E})}
\newcommand{\Dir}{(\mathcal{E},D(\mathcal{E}))}

\title{Construction of $\mathcal L^p$-strong Feller Processes via Dirichlet Forms and Applications to Elliptic Diffusions}
\author{Benedict Baur, Martin Grothaus, Patrik Stilgenbauer}

\begin{document}

\maketitle

\begin{abstract}
We provide a general construction scheme for $\mathcal  L^p$-strong Feller processes on locally compact separable metric spaces. Starting from a regular Dirichlet form and specified regularity assumptions, we construct an associated semigroup and resolvents of kernels having the $\mathcal L^p$-strong Feller property. They allow us to construct a process which solves the corresponding martingale problem for \emph{all} starting points from a \emph{known} set, namely the set where the regularity assumptions hold.
We apply this result to construct elliptic diffusions having locally Lipschitz matrix coefficients and singular drifts on general open sets with absorption at the boundary. In this application elliptic regularity results imply the desired regularity assumptions.
\end{abstract}

\textbf{\underline{AMS classification (2000):}} Primary: 60J25, 31C25; Secondary: 60J60\\

Diffusions processes

\section{Introduction}
In recent years the theory of Dirichlet forms has proven to be a powerful tool for the construction and analysis of stochastic processes going beyond the classical Feller theory. However, the general construction theory for processes from Dirichlet forms yields solutions to the martingale problem outside an exceptional set of starting points only. In general, this set cannot be explicitly specified and in particular need not to be empty. In previous works it turned out that using some elliptic regularity results one obtains a corresponding semigroup of transition kernels and resolvent of kernels possessing enough regularity to construct an associated process, which solves the martingale problem for starting points from an explicitly known set. In \cite{AKR03} distorted Brownian motion on $\R^d$, $d \in \N$, is constructed using elliptic regularity results from \cite{BKR97}. There the admissible starting points are those, where the drift is not singular. In \cite{FG07} this is generalized to the construction of distorted Brownian motion with reflection on domains having smooth boundary except for a known set of capacity zero. There one has to exclude the non-smooth boundary points.\\

The aim of this article is to provide a general construction result for $\mathcal L^p$-strong Feller processes and to apply this to construct elliptic diffusions with singular drifts and locally Lipschitz matrix coefficients.
In the first part of this article we present the construction result for an $\mathcal L^p$-strong Feller process on locally compact metric spaces. The starting point is a symmetric, regular, strongly local Dirichlet form. Additionally, we assume certain regularity on functions in the domain of the associated $L^p$-generator $(L_p,D(L_p))$ for some $p > 1$, see Condition \ref{CondRegularityDomain} below.

The construction result is proven in Section \ref{SectionConstructionKernel} and Section \ref{SectionProcess}. It is based on techniques developed in \cite{AKR03} and \cite{Doh05}. The construction scheme behind is similar to the construction of classical Feller processes. The process then solves the martingale problem for $(L_p,D(L_p))$ for every starting point from a known set.\\

In the second part we apply the results from the first part to construct elliptic diffusions with locally Lipschitz continuous elliptic  matrices and singular drifts on general open sets $\Om \subset \R^d$, $d \in \N$, with absorption at the boundary. For $\Om = \R^d$ and the identity matrix this reproduces the result of \cite{AKR03}. For smooth elliptic matrices first results have been obtained in \cite{Hen08}. Here, however, we consider the case of less smooth matrices. In Section \ref{Dirichlet_form} we construct a gradient Dirichlet form on $\Om$. In Section \ref{regularity_results} we provide an elliptic regularity result based on \cite{BKR01}. This supplies the desired regularity assumptions for our application. In Section \ref{SecConseqRegularity} we apply the results from Section \ref{SectionProcess} to construct elliptic diffusions on $\Om$. The allowed starting points are those, where the drift is not singular.\\
\newline
Let us now state the main results of this article. We refer for notation and results on Dirichlet forms to \cite{FOT94} and \cite{MR92}.
Throughout Section \ref{SectionConstructionKernel} and \ref{SectionProcess} we fix a metric space $(E,d)$, a Borel measure $\mu$ on the Borel $\sigma$-algebra $\mathcal B(E)$ and a symmetric Dirichlet form $(\mathcal E,D(\mathcal E))$. We assume the following conditions.

\begin{condition} \label{CondMetricSpace}
$ $
\begin{enumerate}[label=(\roman{*}),ref=(\roman{*}]
\item $(E,d)$ is a locally compact separable metric space.
\item $\mu$ is locally finite with full topological support.
\item $(\mathcal E,D(\mathcal E))$ is regular and strongly local.
\end{enumerate}
\end{condition}
By the Beurling-Deny theorem there exists an associated strongly continuous contraction semigroup on $L^r(E,\mu)$ ($L^r$-s.c.c.s) $(T^r_t)_{t > 0}$ with generator $(L_r,D(L_r))$ for every $1 \le r < \infty$, see \cite[Prop.~1.8]{LS96} and \cite[Rem.~1.3]{LS96}. If $r>1$ then $(T^r_t)_{t > 0}$ is the restriction of an analytic semigroup, see \cite[Rem.~1.2]{LS96}. Here associated means that for $f \in L^1(E,\mu) \cap L^\infty(E,\mu)$ it holds $T^2_t f = T^r_t f$ for every $t \ge 0$, where $(T^2_t)_{t \ge 0}$ is the unique $L^2$-s.c.c.s associated to $(\mathcal E,D(\mathcal E))$.
\medskip

Additionally we assume the following conditions.
\begin{condition} \label{CondRegularityDomain}
$ $\\
There exists a Borel set $E_1 \subset E$ with $\text{cap}_{\mathcal E}(E \setminus E_1)=0$ and $p > 1$ such that
\begin{enumerate}[label=(\roman{*}),ref=(\roman{*}] 
\item $D(L_p) \emb C(E_1)$ and the embedding is locally continuous, i.e.,~for $x \in E_1$ there exists an $E_1$-neighborhood $U$ and a constant $C_1 = C_1(U) < \infty$ such that
\begin{align}
\sup_{y \in U} |u(y)| \le C_1 \Vert u \Vert_{D(L_p)} \quad \text{for all} \, u \in D(L_p). \label{EqEmbeddingCont}
\end{align}
\item For each point $x \in E_1$ there exists a sequence of functions $(u_n)_{n \in \N}$ in $D(L_p)$ such that
\begin{enumerate}
 \item Either $\{ u^2_n \, | \, n \in \N \} \subset D(L_p)$ or $0 \le u_n \le 1$ and $u_n(x)=1$ for all $n \in \N$.
 \item The sequence $(u_n)_{n \in \N}$ is point separating in $x$.
\end{enumerate}
\end{enumerate}
\end{condition}
Here $C(S)$ denotes the space of all continuous functions on a topological space $S$. By $\Vert \cdot \Vert_{D(L_p)}$ we denote the graph norm of $(L_p,D(L_p))$. Point separating in $x$ means, that for every $y \neq x$ there exists $u_n$ such that $u_n(y)=0$ and $u_n(x)=1$.
Denote by $E^\Delta := E \cup \{ \Delta \}$ the one-point compactification of $E$, endowed with the Alexandrov topology. For $u \in D(L_p)$ we denote the continuous version from Condition \ref{CondRegularityDomain} by $\widetilde u$.
Under Condition \ref{CondMetricSpace} and Condition \ref{CondRegularityDomain} we obtain the following theorem.
\begin{theorem} \label{TheoProcess}
There exists a diffusion process (i.e.,~a strong Markov process having continuous sample paths) $\mathbf M = (\mathbf \Om,\mathcal F, (\mathcal F_t)_{t \ge 0}, (\mathbf X_t)_{t \ge  0},(\mathbf P_x)_{x \in E_1 \cup \{ \Delta \}})$ with state space $E_1$ and cemetery $\Delta$, the Alexandrov point of $E$. The transition semigroup $(P_t)_{t \ge 0}$ is associated to $(T^2_t)_{t \ge 0}$ and is $\mathcal L^p$-strong Feller, i.e.,~$P_t \mathcal L^p(E,\mu) \subset C(E_1)$ for $t>0$. The process has continuous paths on $[0,\infty)$ and it solves the martingale problem associated to $(L_p,D(L_p))$, i.e.,
\begin{align*}
M_t^{[u]}:=\widetilde u({\mathbf X}_t) - \widetilde u(x) - \int_0^t{L_p u({\mathbf X}_s)~ds},~t\geq0,
\end{align*}
is an $(\mathcal{F}_t)$-martingale under $\mathbf{P}_x$ for all $u \in D(L_p)$, $x \in E_1$.
\end{theorem}
Here $(P_t)_{t \ge 0}$ being associated to $(T^2_t)_{t \ge 0}$ means that $P_t f$ is a $\mu$-version of $T^2_t f$ for $f \in \mathcal L^1(E,\mu) \cap \mathcal B_b(E)$ (the space of Borel-measurable bounded functions). With $\mathcal L^p(E,\mu)$ we denote the space of all $p$-integrable functions on $(E,\mu)$.
\begin{remark}
The continuity holds with respect to the Alexandrov topology of the one-point compactification of $E$ to $E^\Delta$. This means that the process has continuous paths in $E$ and reaches $\Delta$ only by leaving continuously every compact set of $E$. 
\end{remark}
The theorem is proven in Section \ref{SectionProcess}, see page \pageref{proofTheoDiffProcess} and Theorem \ref{TheoremMartingaleSolution} below. Further useful properties of the constructed process are proven in Theorem \ref{TheoPropHit} below.

Under additional conditions, the corresponding resolvent of  kernels $(R_\lambda)_{\lambda > 0}$ are even \textit{strong Feller}, i.e., $R_\lambda \mathcal B_b(E) \subset C(E_1)$. More precisely, we have the following theorem. 
\begin{theorem} \label{ClassicalStrongFeller}
Assume the following conditions.
\begin{enumerate}[label=(\roman{*}),ref=(\roman{*}]
 \item For every $x \in E_1$ there exists a neighborhood $U \subset E_1$ such that for the closure in $E$ it holds $\overline{U} \subset E_1$ and $\overline{U}$ is compact.
 \item For every sequence $(u_n)_{n \in \N}$ in $D(L_p)$ such that $((1 - L)u_n)_{n \in \N}$ is uniformly bounded in the $\Vert \cdot \Vert_{L^\infty}$-norm it holds that $(u_n)_{n \in \N}$ is equicontinuous.
\end{enumerate}
Then $(R_\lambda)_{\lambda > 0}$ is strong Feller.
If additionally $(\mathcal E,D(\mathcal E))$ is conservative, then $\mathbf M$ from Theorem \ref{TheoProcess} is conservative.
\end{theorem}
For the proof see Section \ref{SectionConstructionKernel} (page \pageref{ProofClassicalStrongFeller}).
\begin{remark}
In \cite{AKR03} it is shown, that the strong Feller property of $(R_\lambda)_{\lambda > 0}$ and conservativity of $(\mathcal E,D(\mathcal E))$ imply that $(P_t)_{t > 0}$ is strong Feller. The proof generalizes to the case considered here.
\end{remark}

In Section \ref{SecConseqRegularity} we apply this result to construct elliptic diffusions. For this we fix an open subset $\Om \subset \R^d$, $d \ge 2$, a matrix-valued mapping $A = (a_{ij})_{i,j=1}^d : \Om \to \R^{d\times d}$ a Borel measure $\mu$ with density $\varrho : \Om \to \R^+$, i.e., $\mu := \varrho dx$. Furthermore, we fix $p \in \N$. 
We consider the following pre-Dirichlet form
\begin{align}
&\E(f,g) := \int_{\Omega} ~ \big(A \nabla f,\nabla g\big)_{\text{euc}} ~ d\mu = \int_{\Omega} ~ \sum_{i,j=1}^d a_{ij} \, \pdi f \, \pdj g ~ d\mu,~ f,g \in \mathcal D:=C_c^\infty(\Omega), \label{EqDirichletFormIntroduction}
\end{align}
where $C_c^\infty(\Omega)$ denotes the space of smooth functions on $\Om$ with compact support in $\Om$.

We assume the following conditions.
\begin{condition} \label{AssEllipConstruction}
Assume $p > d$.
\begin{enumerate}[label=(\roman{*}),ref=(\roman{*}]
\item $\sqrt{\varrho} \in H^{1,2}_{\text{loc}}(\Omega), ~ \varrho > 0 ~ dx$-a.e.
\item $\frac{\nabla\varrho}{\varrho} \in L^p_{\text{loc}}(\Omega,\mu)$.
\item $A$ is measurable, symmetric and locally strictly elliptic on $\Omega$ $dx$-a.e. Furthermore, $a_{ij} \in H^{1,\infty}_{loc}(\Omega)$, $1 \leq i,j \leq d$.
\end{enumerate}
\end{condition}
Here $H^{m,p}(\Om) \subset L^p(\Om,dx)$ denotes the Sobolev space of $m$-times weakly differentiable functions with $L^p(\Om,dx)$-integrable derivatives.
The notation $f \in M_{\text{loc}}(U)$, $U \subseteq \mathbb{R}^d$ open, $M$ some Sobolev or Lebesgue space, is understood in the sense that $f$ is measurable and locally an element of $M$. Then we obtain the following proposition.
\begin{proposition} \label{PropDirichletform}
The density $\varrho$ has a unique continuous $\mu$-version.
$(\E,D)$ is closable on $L^2(\Omega,\mu)$ with closure $\Dir$. The closure is a symmetric, strongly local and regular Dirichlet form on $L^2(\Omega,\mu)$.
%Moreover, for each $1 \le r < \infty$ there exists an associated $L^r$-s.c.c.s with generator $(L_r,D(L_r))$.
For $p > d$ as in Condition \ref{AssEllipConstruction} it holds $C^\infty_c(\Om) \subset D(L_p)$.
Furthermore, we have $\text{cap}_{\mathcal E}(\{ \varrho = 0 \}) = 0$.
\end{proposition}
This is proven in Section \ref{Dirichlet_form}, see Remark \ref{RmHolderContinuousVersionRho}, Proposition \ref{Pp Dirichlet form}, Proposition \ref{Pp Generator} and Proposition \ref{Pp Capacity} below.
\begin{remark}
We fix the continuous version of $\varrho$ provided by Proposition \ref{PropDirichletform} and denote it also by $\varrho$.
\end{remark}

Based on the elliptic regularity results from Section \ref{regularity_results} we obtain the desired regularity for functions in $D(L_p)$, see Section \ref{SecConseqRegularity} below.
\begin{theorem} \label{Th Resolvent Continuous}
For $p$ as in Condition \ref{AssEllipConstruction} it holds $D(L_p) \subset C(\{ \varrho > 0\})$, in the sense that for each $u \in D(L_p)$ there exists a unique $\mu$-version on $\{ \varrho > 0 \}$ which is continuous on $\{ \varrho > 0\}$. This version is denoted by $\widetilde u$.
For each $x \in \{ \varrho > 0\}$ and $r > 0$ such that $B := B_r(x)$ (Ball of radius $r$ around $x$) has closure in $\{ \varrho > 0\}$  
it holds $\widetilde u \in C^{0,\beta}(B)$ with $\beta := 1 - \frac{d}{p}$ and
there exists $K < \infty$ independent of $u$ such that 
\begin{align}
\Vert \widetilde{u} \Vert_{C^{0,\beta}(\overline{B})} \le K(\Vert u \Vert_{L^p(B,\mu)} + \Vert L u \Vert_{L^p(B,\mu)}) \le K \Vert u \Vert_{D(L_p)}. \label{Estimate Th Resolvent Continuous}
\end{align}
\end{theorem}
Here $C^{0,\beta}(\overline{B})$ denotes the space of H"older continuous functions with index $\beta$ on $\overline{B}$. We denote by $|f|_{C^{0,\beta}(\overline{B})}$ the H"older coefficient and by  $\Vert f \Vert_{C(\overline{B})}$ the sup norm of $f$ in $\overline{B}$. Then we define $\Vert f \Vert_{C^{0,\beta}(\overline{B})} := \Vert f \Vert_{C(\overline{B})} + |f|_{C^{0,\beta}(\overline{B})}$.
\\
Finally, using Theorem \ref{TheoProcess} with $E_1 = \{ \varrho > 0 \} \subset \Om$ we prove the following theorem in Section \ref{SecConseqRegularity}.  
\begin{theorem} \label{DiffProcess}
There exists a diffusion process
\begin{align*}
\mathbf{M} = \left( \mathbf{\Omega}, ~\mathcal {F}, (\mathcal{F}_t)_{t\geq 0},({\mathbf X}_t)_{t\geq 0}, (\mathbf{P}_x)_{x \in \{ \varrho > 0 \} \cup \{ \Delta \} }\right) 
\end{align*}
with state space $\{ \varrho > 0 \}$ and cemetery $\Delta$, the Alexandrov point of $\Om$. The transition semigroup $(P_t)_{t > 0}$ 
is associated to $(T^2_t)_{t > 0}$ and is $\mathcal L^p$-strong Feller, i.e.,~$P_t \mathcal L^p \subset C(\{ \varrho > 0 \})$. Furthermore, the associated resolvent of kernels $(R_\lambda)_{\lambda > 0}$ are strong Feller, i.e.,~$R_\lambda \mathcal B_b \subset C(E_1)$.
The process has continuous paths on $[0,\infty)$ and it yields a solution to the martingale problem for $(L_p,D(L_p))$, i.e.,
\begin{align*}
M_t^{[u]}:=\widetilde u({\mathbf X}_t) - \widetilde u(x) - \int_0^t{L_p u({\mathbf X}_s)~ds},~t\geq0,
\end{align*}
is an $(\mathcal{F}_t)$-martingale under $\mathbf{P}_x$ for all $u \in D(L_p)$, $x \in \{ \varrho > 0 \} \cup \{ \Delta \}$. In particular, this holds for $u \in C^\infty_c(\Om)$.
\end{theorem}
Here $\widetilde u$ denotes the continuous version on $E_1$ of $u$.

\section{Construction of $\mathcal L^p$-strong Feller kernels} \label{SectionConstructionKernel}
%Let $(E,d)$ be a locally compact separable metric space. Let $\mu$ be a measure on $(E,\mathcal B(E))$, which is finite on compact sets and $\supp(\mu) = E$, i.e. every non-empty open set has strictly positive measure.
%Furthermore, assume that there exists a strongly local regular symmetric Dirichlet form $(\mathcal E,D(\mathcal E))$ on $L^2(E,\mu)$ and a measurable set $E_1$ with $cap_{\mathcal E}(E \setminus E_1)=0$ having the following properties: There exists $p \in \N$ such that for the $L^p$-s.c.c.s $(T^p_t)_{t \ge 0}$ associated to $\mathcal E$ it holds that $D(L_p) \emb C(E_1)$ and the embedding is locally continuous, that is for each $x \in E_1$ there exists a neighborhood $U$ of $x$ w.r.t the topology on $E$ with $U \subset E_1$ and a constant $C < \infty$ such that
%\begin{align}
%\sup_{x \in U} |u(x)| \le C \Vert u \Vert_{D(L_p)}. \label{EqEmbeddingCont}
%\end{align}
%Here $\Vert \cdot \Vert_{D(L_p)}$ denotes the graph-norm of the $L^p$-generator.
We start with the construction of a semigroup of kernels $(P_t)_{t > 0}$ and resolvent of kernels $(R_\lambda)_{\lambda > 0}$, which yield a $\mu$-version of $(T^p_t)_{t > 0}$ and $(G^p_\lambda)_{\lambda > 0}$. For this we assume Condition \ref{CondMetricSpace} and Condition \ref{CondRegularityDomain} and fix a $p>1$ as in Condition \ref{CondRegularityDomain}. The concepts and proofs from \cite{AKR03} are generalized to the abstract setting. 
\begin{remark}
The restriction of $\mu$ to $\mathcal B(E_1)$ is also strictly positive on non-empty open sets. Indeed, let $\widetilde U \subset E_1$ be non-empty and open w.r.t.~the trace topology. Then there exists $U \subset E$ open, such that $\widetilde U = U \cap E_1$. So $U \setminus \widetilde U \subset E \setminus E_1$, the latter is of capacity zero and hence has also $\mu$ measure equal to zero. So $\mu(\widetilde U) = \mu(U) > 0$.  
In particular, if $\tilde u$ is continuous on $E_1$ and equal to zero $\mu$-almost everywhere in $E_1$, then $\tilde u$ is equal to zero on $E_1$.
This implies, that if $u \subset \mathcal B(E)$ has a continuous version on $E_1$, then this version is unique on $E_1$.
\end{remark}
For the associated $L^p$-resolvent $(G^p_\lambda)_{\lambda > 0}$ it holds $D(L_p) = G^p_\lambda L^p(E,\mu)$. So for $f \in L^p(E,\mu)$, $G^p_\lambda f$ has a unique continuous version denoted by $\widetilde{G^p_\lambda f}$. The boundedness of $G_\lambda : L^p(E,\mu) \to D(L_p)$ together with $\eqref{EqEmbeddingCont}$ yields a constant $C_2 = C_2(\lambda,U) < \infty$, $U$ as in Condition \ref{CondRegularityDomain}, such that
\begin{align}
 \sup_{y \in U} |\widetilde {G^p_\lambda f} (y) | \le C_1 \Vert G^p_\lambda f \Vert_{D(L_p)} \le C_2 \Vert f \Vert_{L^p(E,\mu)} \quad \text{for all} \ f \in L^p(E,\mu). \label{EqResolventEst}
\end{align}
Since $p>1$ the $L^p$-semigroup $(T^p_t)_{t > 0}$ is the restriction of an analytic semigroup. Thus $T^p_t L^p(E,\mu) \subset D(L_p)$ for $t > 0$. So, for $u \in L^p(E,\mu)$, $t > 0$, $T^p_t f$ has a unique continuous version, denoted by $\widetilde{T^p_t f}$. 
From \cite[Theo.~5.2]{Pa83} we get that there exists a constant $C_3 < \infty$ such that
\begin{align*}
\Vert T^p_t f \Vert_{D(L_p)} \le \left( 1+\frac{C_3}{t} \right) \Vert f \Vert_{L^p(E,\mu)} \quad \text{for all} \, f \in L^p(E,\mu). 
\end{align*}
Together with \eqref{EqEmbeddingCont} it follows then that there exists a constant $C_4 = C_4(t,U) < \infty$, depending on $t$, such that
\begin{align}
 \sup_{y \in U} | \widetilde{T^p_t f}(y) | \le C_1 \Vert \widetilde{T^p_t f} \Vert_{D(L_p)} \le C_4 \Vert f \Vert_{L^p(E,\mu)} \quad \text{for all} \, f \in L^p(E,\mu). \label{EqSemigroupPointwise}
\end{align}
In particular, $L^p$-convergence of a sequence $(f_n)_{n \in \N}$ to $f$ implies pointwise convergence of $\widetilde{T^p_t f_n}(x)$ to $\widetilde{T^p_t f}(x)$, $t > 0$, and pointwise convergence of $\widetilde {G^p_\lambda f_n}(x)$ to $\widetilde {G^p_\lambda f}(x)$, $\lambda > 0$, for $x \in E_1$.

If $u \in D(L_p)$ then $\Vert T_t u \Vert_{D(L_p)} \le \Vert u \Vert_{D(L_p)}$ and hence
\begin{align}
 \sup_{y \in U} | \widetilde{T^p_t u}(y) | \le C_1 \Vert u \Vert_{D(L_p)}. \label{EqSemigroupGraphnorm}
\end{align}

Most of the time we omit the upper index $p$ when applying the $L^p$-semigroup or resolvent to a function $f \in L^p(E,\mu)$.

The following well-known lemma is useful for monotone class arguments.
% Kurze Version, ausf"uhrliche in anderem File
\begin{lemma} \label{LemmaSigmaAlgebraOpenFinite}
Let $(E,d)$ be a locally compact separable metric space. Then there exists a sequence of compact sets $(K_n)_{n \in \N}$, $K_n \subset \interior{K}_{n+1}$ (the interior of $K_{n+1}$), with $E = \bigcup_{n \in \N} K_n$. Furthermore, if $\mu$ is a locally finite measure on $\mathcal B(E)$, then $\mathcal B(E)$ is generated by the open sets of finite measure.
\end{lemma}
%\begin{proof}
%Using the fact the $(E,d)$ is separable and locally compact one constructs a countable covering of balls with compact closure. The usual vast-out argument for open sets yields the existing of compact sets $K_n$ with $K_{n} \subset \interior{K_{n+1}}$, $n \in \N$, covering $E$. In particular also $\interior{K_{n}}$ cover $E$ and have finite measure since $\mu(K_{n+1}) < \infty$.
%For any open set $U$, $U \cap \interior{K_n}$ is also open and has finite measure. So $U$ is contained in the $\sigma$-algebra generated by the open sets of finite measure. 
%\end{proof}
For the proof of the existence of such a covering, see e.g. \cite[Cor.~2.77]{CB06}. The second statement follows then directly.

For many constructions we need the Functional Monotone Class Theorem, see e.g. \cite[Ch.~0, Theo.~2.3]{BG68}.
\begin{theorem} \label{TheoFMCT}
Let $(E,\mathcal B)$ be a measurable space and $\mathcal L \subset \mathcal B$ be an intersection-stable generator of $\mathcal B$.
Let $\mathcal H \subset \mathcal B_b(E)$ be a subset having the following three properties:
\begin{enumerate}[label=(\roman{*}),ref=(\roman{*}]
\item $\mathcal H$ is a vector space over $\R$. \label{VS}
\item $\mathcal H$ contains $1_E$ as well as the indicator function $1_F$ for $F \in \mathcal L$. \label{SG}
\item If $f_n \in \mathcal H$, $f_n \ge 0$, $n \in \N$, such that $f_n \uparrow f$ as $n \to \infty$ and $f$ is bounded, then $f \in \mathcal H$. \label{MC}
\end{enumerate}
Then $\mathcal H = \mathcal B_b(E)$.
\end{theorem}
%\begin{proof}
%First set $\mathcal F := \{ F \in \mathcal B \, | \, 1_F \in \mathcal U \}$. We prove that this is a monotone class (also called a Dynkin-system). First note that by assumption (SG) $E \in \mathcal F$. Moreover, for $F_1,F_2 \subset \mathcal F$ with $F_1 \subset F_2$ we have $1_{F_2 \setminus F_1} = 1_{F_2} - 1_{F_1}$ so $1_{F_2 \setminus F_1} \in \mathcal U$ by (VS). If $F_n \subset \mathcal F$, $n \in \N$, are pairwise disjoint then for $F := \bigcup_{n \in \N} F_n$ it holds that $\sum_{n=1}^N 1_{F_n} \uparrow 1_F$. So, by (MC) it holds $1_F \in \mathcal U$. So, $\mathcal F$ is a monotone class, which contains by (SG) an intersection-stable generator of $\mathcal B$. Thus by the Monotone Class Theorem it follows that $\mathcal F$ contains $\mathcal B$. Since every bounded positive measurable function can be approximated monotonically by linear combinations of indicator functions the claim follows using (VS) and (MC).
%\end{proof}

\begin{corollary} \label{CoroTheoFMCT}
Assume that $(E,d)$ is a locally compact separable metric space, $\mu$ a locally finite measure on $\mathcal B(E)$. Then the conclusion of Theorem \ref{TheoFMCT} remains true, if instead of (ii) we assume either \\

(ii') \ For $U$ open, $\mu(U) < \infty$, it holds $1_U \in \mathcal H$, \label{SG'} or \\

(ii'') \ For $U$ open, $\mu(U) < \infty$, there exist $f_n \in \mathcal H$ such that $f_n \uparrow 1_U$. \label{SGA}\\

In particular, if $C_b(E) \cap L^p(E,\mu) \subset \mathcal H$ and $\mathcal H$ has \ref{TheoFMCT}(i) and \ref{TheoFMCT}(ii), then $\mathcal H = \mathcal B(E)$.
\end{corollary}
\begin{proof}
(ii') implies \ref{TheoFMCT}(ii): We have $E = \bigcup_{n \in \N} \interior{K}_n$ with $K_n$ as in Lemma \ref{LemmaSigmaAlgebraOpenFinite}. Thus, the open sets of finite measure are an intersection stable-generator. Hence, we can choose as $\mathcal L$ the system of all open sets with finite measure. Moreover, $1_{\interior{K_n}} \uparrow 1_E$. So, together with \ref{TheoFMCT}(iii) it follows $1_E \in \mathcal H$. Altogether, (ii') implies \ref{TheoFMCT}(ii).\\
(ii'') implies (ii'): Follows directly by \ref{TheoFMCT}(iii).\\
For the proof of the last claim, recall that for every open set $U$, there exists a sequence of continuous functions $(f_n)_{n \in \N}$ with $f_n \uparrow 1_U$. So $f_n \in C_b(E)$ and if $\mu(U) < \infty$ then also $f_n \in L^p(E,\mu)$. So, if $C_b(E) \cap L^p(E,\mu) \subset \mathcal H$, then $\mathcal H$ fulfills (ii'').
\end{proof}

Now we construct the semigroup of kernels $(P_t)_{t \ge 0}$ on $E_1 \times \mathcal B(E)$. For this we need the following Lemma.
\begin{lemma} \label{LemmaDaniellIntegral}
Let $t>0$, $x \in E_1$. Then the map
\begin{align}
\mathcal L^p(E,\mu) \ni \,  f  \mapsto \widetilde{T_t f}(x)  \in \R,  \label{EqMapTildeT}
\end{align}
is a Daniell-integral, cf.~\cite[Def.~39.1]{Bau78}, and there exists a unique positive measure $P_t(x,dy)$ on $\mathcal B(E)$ such that
\begin{align}
 \widetilde{T_t f}(x) = \int_{E} f(y) P_t(x,dy). \label{EqDaniellMap}
\end{align}
\end{lemma}
\begin{proof}
By positivity of $T_t f$ and continuity on $E_1$ we have that $\widetilde{T_t f}(x) \ge 0$ for every $x \in E_1$, if $f \ge 0$. Using linearity of $T_t$ and continuity of $\widetilde{ T_t f}$ we get that the mapping \eqref{EqMapTildeT} is also linear. Moreover, if $f_n \downarrow 0$, $\mu$-a.e. this convergence also holds in $L^p(E,\mu)$ by Lebesgue's dominated convergence. So $\eqref{EqSemigroupPointwise}$ yields $\widetilde{T_t f}(x) \to 0$. By positivity this convergence is also monotone. Thus, the map is a Daniell-integral. By \cite[Satz 39.4]{Bau78} there exists a positive measure denoted by $P_t(x,dy)$ on $\sigma(\mathcal L^p)$ (the $\sigma$-algebra generated by $\mathcal L^p$) such that $\eqref{EqDaniellMap}$ holds.
Note that for every set $M \in \mathcal B(E)$ of finite measure it holds $1_M \in \mathcal L^p(E,\mu)$. By Lemma \ref{LemmaSigmaAlgebraOpenFinite} $\mathcal B(E)$ is generated by the open sets of finite measure. Hence, $\mathcal B(E) = \sigma(\mathcal L^p(E,\mu))$.
Moreover, the measure is unique, since the open sets of finite measure are an intersection-stable generator.
\end{proof}

\begin{remark}
Note that the map $\eqref{EqMapTildeT}$ is formulated on the $\mathcal L^p$-functions rather than on the $\mu$-equivalence classes from $L^p(E,\mu)$. However, the operator $T_t$ respects $\mu$-equivalence classes. So, two different representatives of an element in $L^p(E,\mu)$ lead to the same equivalence class $T_t f$ and to the same unique continuous version $\widetilde {T_t f}$. So the mapping $\eqref{EqMapTildeT}$ is also well-defined as a mapping $L^p(E,\mu) \to \R$.
\end{remark}

Now we prove some facts about the kernels $P_t(x,dy)$, $t \ge 0$, $x \in E_1$.

\begin{definition} \label{DefKernelPt}
For $t > 0$, $x \in E_1$ and $f \in \mathcal L^1(E, P_t(x,dy)) \cup \mathcal B^+(E)$ we define
\begin{align*}
 P_t f (x) := \int_{E} f(y) P_t(x, dy)
\end{align*}
and $P_0 f(x) := f(x)$.
\end{definition}

\begin{theorem} \label{TheoremKernelSemigroup}
$ $
\begin{enumerate}
\item[(i)] 
Let $x \in E_1$.
It holds $P_t 1_E (x)\le 1$.
There exists a $\mathcal B^*(E_1 \times E)$-measurable map $(x,y) \mapsto p_t(x,y)$, $t > 0$ such that $P_t(x,dy) = p_t(x,y) \mu(dy)$. Furthermore, $P_t(x,E \setminus E_1) = 0$, so the kernels defined in Definition \ref{DefKernelPt} can be considered as kernels on $E_1$, denoted below by the same symbol.
Moreover, $\mathcal L^p(E,\mu) \subset \mathcal L^1(E,P_t(x,dy))$ and $P_t f (x)= \widetilde{T^p_t f} (x)$ for all $x \in E_1$, $f \in \mathcal L^p(E,\mu)$, i.e. $P_t f$ is the unique continuous version of $T^p_t f$. 
% The integral $P_t f$ coincides for all $\mu$-version of $f$, so $P_t$ is also well-defined as a map $L^p(E;d\mu) \to C(E_1)$.
\item[(ii)] $(P_t)_{t > 0}$ is a semigroup of kernels on $E_1$ which is $\mathcal L^p$-strong Feller, i.e. $P_t f \in C(E_1)$ for all $t>0$, $f \in \mathcal L^p(E, \mu)$. Moreover, $P_{t+s} f = P_t P_s f$ for $f \in \mathcal L^p(E,\mu)$.
\item[(iii)] For $f \in L^p(E,\mu)$ and $s > 0$
\begin{align*}
\underset{t \to 0}{\lim} \, P_{t+s} f (x) = P_s f(x) \quad \text{for all} \  x \in E_1.
\end{align*}
For $s=0$ this holds for all $f \in D(L_p)$.
\item[(iv)] $(P_t)_{t > 0}$ is a measurable semigroup on $E_1$, i.e. for $f \in \mathcal B_b(E)$ the map $(t,x) \mapsto P_t f(x)$ is $\mathcal B([0,\infty) \times E_1)$-measurable. This holds also for $f \in \mathcal L^p(E,\mu)$.
\end{enumerate}
\end{theorem}
Here $\mathcal B^*(E_1 \times E)$ denotes the completion of the Borel $\sigma$-algebra $\mathcal B(E_1 \times E)$ with respect to the product measure $dx \otimes \mu$.
\begin{proof}
(i): Let $F_k \subset E$, $k \in \N$, be sequence with $\mu(F_k) < \infty$ and $E = \bigcup_{k \in \N} F_k$. Then $T^p_t 1_{F_k} \le 1$ $\mu$-a.e on $E$. By continuity we have $P_t 1_{F_k}(x) = \widetilde{T^p_t 1}_{F_k} (x) \le 1$ for every $x \in E_1$. By monotone convergence it holds $P_t 1_E(x) = \sup_{k \in \N} P_t 1_{F_k}(x) \le 1$ for $x \in E_1$ and $t>0$.
Let $N \in \mathcal B(E)$ with $\mu(N) = 0$. Then $1_N = 0$ $\mu$-a.e., so $T^p_t 1_N = 0$ $\mu$-a.e. So by continuity $P_t(x,N) = P_t 1_N (x)= \widetilde{T^p_t 1}_N (x)= 0$ for every $x \in E_1$. Hence, $P_t(x,dy)$ is absolutely continuous w.r.t.~$\mu$. Thus, by the Radon-Nikodym theorem there exists a map $p_t(x,y)$ such that $p_t(x,y)$ is the density of $P_t(x,dy)$ for all $x \in E_1$ and $t>0$. 
To prove measurability note that for $U \in \mathcal B(E)$ with $\mu(U) < \infty$ the mapping $x \to P_t(x,U) = \widetilde{T^p_t 1}_U(x)$ is measurable in $x$, hence by monotone approximation using Lemma \ref{LemmaSigmaAlgebraOpenFinite} also for arbitrary $U \in \mathcal B(E)$. The existence of a $\mathcal B^*(E_1 \times E)$-measurable version of $p_t(x,y)$ follows by constructing the density as the Radon-Nikodym derivative along a partition of $E$, see \cite[Theo.~2.5]{Doo53} and \cite[Exa.~2.7]{Doo53}.
Since $\text{cap}_{\mathcal E}(E \setminus E_1) = 0$ we have $\mu(E \setminus E_1)=0$. So $P_t(x,E \setminus E_1) = 0$. 
Now let $f \in \mathcal L^p(E,\mu)$ and consider $f^+$. Then for all $x \in E_1$, $t>0$ we have by construction of $P_t$
\begin{align*}
P_t f^+ (x) = \widetilde{T^p_t f^+} (x) < \infty.
\end{align*}
Thus $f^+ \in \mathcal L^1(E,P_t(x,dy))$. The same reasoning works for $f^-$, thus $f \in \mathcal L^1(E,P_t(x,dy))$ and $P_t f (x) = \widetilde{T^p_t f} (x)$. \newline
%By construction $P_t f$ is independent of the $\mu$-version of $f$, so $P_t : L^p(E,d\mu) \to C(E_1)$ is well-defined.
(ii): Let $u \in \mathcal L^p(E,\mu)$. Then by the $L^p(E,\mu)$-semigroup property and continuity we have for all $x \in E_1$, $t,s \ge 0$
\begin{align*}
\widetilde{T_{t+s} u} (x) = \widetilde{T_t (T_s u)} (x).
\end{align*}
Thus,
\begin{align}
P_{t+s} u (x) = P_t (P_s u) (x) \quad \text{for all} \ x \in E_1. \label{EqSemigroupPropertyProof}
\end{align}
Note that this holds in particular for $u = 1_U$, $U$ open with $\mu(U) < \infty$. To prove it for Borel bounded functions observe that the system of functions in $\mathcal B_b(E)$ for which property $\eqref{EqSemigroupPropertyProof}$ holds is a vector space satisfying \ref{TheoFMCT}(iii). So by Corollary \ref{CoroTheoFMCT} it follows that $\eqref{EqSemigroupPropertyProof}$ holds also for $u \in \mathcal B_b(E)$. The other statements are clear by construction.
\newline
(iii): First, we prove the statement for $s=0$ and $u \in D(L_p)$.
Note that the $L^p$ semigroup $(T_t)_{t \ge 0}$ is strongly continuous on $D(L_p)$ w.r.t.~the graph norm of $(L_p,D(L_p))$.
So for $u \in D(L_p)$ we get using $\eqref{EqEmbeddingCont}$ for $x \in E_1$
\begin{align*}
 |\widetilde {T_t u}(x) - \widetilde{u}(x)| \le C_1 \Vert T_t u -u \Vert_{D(L_p)} \overset{t \to 0}{\longrightarrow} 0.
\end{align*}

Since $P_t u (x) = \widetilde{T^p_t u} (x)$ we get 
\begin{align}
\underset{t \to 0}{\lim} \ P_t u (x) = \underset{t \to 0}{\lim} \ \widetilde{T^p_t u} (x) = u(x) = P_0 u (x) \quad \text{for all} \ x \in E_1. \label{EqKernelContinuity} 
\end{align}
For $s>0$, $f \in L^p(E,\mu)$, we have by analyticity of $(T_t)_{t > 0}$ that $P_s f = \widetilde{T_s f} \in D(L_p)$. By the semigroup property of the kernels $(P_t)_{t \ge 0}$ we have
\begin{align*}
\underset{t \to 0}{\lim} \ P_{t+s} f (x) = \underset{t \to 0}{\lim} \ P_t (P_s f) (x) = P_0 (P_s f) (x) = P_s f (x). 
\end{align*}
(iv): First let $f \in \mathcal L^p(E,\mu)$. Define for $n \in \N_0$, $S_n := \{ k 2^{-n} \, | \, k \in \N_0 \}$, $s^{(n)}_k := k 2^{-n}$, $k \in \N_0$, $M^{(n)}_k = (s^{(n)}_{k-1},s^{(n)}_k]$, $k \in \N$, and $M^{(n)}_0 = \{ 0 \}$. For $t>0$ define $t_n := \min \{ s \in S_n \, | \, t \le s \}$. Clearly $t_n \downarrow t$ as $n \to \infty$.
We define for $n \in \N$
\begin{align*}
P^n f : (t,x) \mapsto P^n_t f(x) := P_{t_n} f(x).
\end{align*}
Then for $A \in \mathcal B(\R)$, it holds
\begin{align*}
(P^n_{\cdot} f(\cdot))^{-1}(A) = \bigcup_{k \in \N_0} M_k^{(n)} \times (P_{s^{(n)}_k} f)^{-1}(A) \in \mathcal B(\R_0^+ \times E_1).
\end{align*}
Thus $P^n f$ is measurable. Now note that for $t=0$ we have $P^n_t f(x) = P_0 f(x)$ and for $t>0$ we have $P^n_t f(x) = P_{t_n} f(x) \overset{n \to \infty}{\longrightarrow} P_t f(x)$ for $x \in E_1$ by $(iii)$.
So $P_t f$ is measurable for $f \in \mathcal L^p(E,\mu)$. Then measurability for general $u \in \mathcal B_b(E)$ follows as in $(ii)$ using Corollary \ref{CoroTheoFMCT}.
\end{proof}
For the resolvent similar statements hold.
\begin{lemma} \label{LemmaDaniellIntegralResolvent}
Let $0 < \lambda < \infty$ and $x \in E_1$. Then the map
\begin{align*}
\mathcal L^p(E,\mu) \ni \, f \mapsto \widetilde{G^p_\lambda f}(x) \in \R
\end{align*}
on $\mathcal L^p(E,\mu)$ is a Daniell-integral, hence there exists a unique positive measure $R_\lambda(x,dy)$ on $\mathcal B(E)$ such that
\begin{align*}
\widetilde{G^p_\lambda f} (x) = \int_{E} f(y) R_\lambda(x,dy) \quad \text{for all} \, \, f \in \mathcal L^p(E,\mu).
\end{align*}
\end{lemma}

\begin{definition} \label{DefKernelR}
For $\lambda > 0$, $x \in E_1$ and $f \in L^1(E, R_\lambda(x,dy)) \cup \mathcal B^+(E)$ we define
\begin{align*}
R_\lambda f (x) := \int_{E} f(y) R_\lambda(x, dy).
\end{align*}
\end{definition}

\begin{theorem} \label{PropKernelResolvent}
$ $
\begin{enumerate}
\item[(i)] For $\lambda > 0$ it holds $\lambda R_\lambda 1 \le 1$. There exists a $\mathcal B^*(E_1 \times E)$-measurable map $(x,y) \mapsto r_\lambda(x,y)$ such that $R_\lambda(x,dy) = r_\lambda(x,y) \mu(dy)$. In particular, $R_\lambda(x,E \setminus E_1) = 0$, so $R_\lambda$ can be considered as kernels on $E_1$, denoted by the same symbol.
Moreover, $\mathcal L^p(E,\mu) \subset \mathcal L^1(E_1,R_\lambda(x,dy))$ and $R_\lambda f (x)= \widetilde{G^p_\lambda f} (x)$ for all $x \in E_1$, i.e. $R_\lambda f$ is the unique continuous version of $\widetilde{G^p_\lambda f}$. In particular, the integral $R_\lambda f$ coincides for all $\mu$-version of $f$.
\item[(ii)] $(R_\lambda)_{\lambda > 0}$ is a resolvent of kernels on $E_1$ which is $\mathcal L^p$-strong Feller, i.e. $R_\lambda f \in C(E_1)$ for all $f \in \mathcal L^p(E,\mu)$. 
% \item[(iii)] $(R_\lambda)_{\lambda > 0}$ is
% Strong Feller, i.e. $R_\lambda f \in C_b(\Omcloplus)$ for $f \in \mathcal B_b(\Omclo)$. 
\item[(iii)] For all $u \in D(L_p)$ and all $x \in E_1$
\begin{align}
R_\lambda u(x) = \int_0^\infty \exp(-\lambda t) P_t u (x) \, dt. \label{EqResolventLaplaceTrafo}
\end{align}
\item[(iv)] For all $u \in D(L_p)$
\begin{align*}
\underset{\lambda \to \infty}{\lim} \lambda R_\lambda u (x) = u(x) \quad \text{for all} \ x \in E_1.
\end{align*}
\end{enumerate}
\end{theorem}
\begin{proof}
The proofs of (i) and (ii) work analogously to those of Theorem \ref{TheoremKernelSemigroup}.\\
(iii): Let $u \in D(L_p)$. Note that by the properties of $G^p_\lambda$ and $T^p_t$ it holds 
\begin{align*}
G_\lambda u (x) = \int_0^\infty \exp(-\lambda t) T_t u \, dt \, (x) = \int_0^\infty \exp(-\lambda t) P_t u \, dt \, (x) \quad \, \text{for} \ \mu\text{-a.e.} \ x.
\end{align*}
The integral is obtained as the $L^p$-limit and hence by dropping to a subsequence also as the $\mu$-a.e.~limit of Riemannian sums. Furthermore, the mapping $[0,\infty) \mapsto P_t u(x)$ is continuous and bounded by $\eqref{EqSemigroupGraphnorm}$, hence the mapping $[0,\infty) \mapsto \exp(-\lambda t) P_t u(x)$ is Lebesgue-integrable and the integral is obtained as the limit of Riemannian sums. Thus, 
\begin{align*}
\int_0^\infty \exp(-\lambda t) P_t u \, dt \, (x) = \int_0^\infty \exp(-\lambda t) P_t u (x) \, dt \quad \, \text{for} \ \mu\text{-a.e.} \ x.
\end{align*}
Thus, we get for almost all $x \in E_1$ 
\begin{align}
R_\lambda u (x) = G_\lambda u(x) = \int_0^\infty \exp(-\lambda t) T_t u \, dt \, (x) = \int_0^\infty \exp(-\lambda t) P_t u (x) \, dt. \label{EqLaplace0}
\end{align}
By $\eqref{EqSemigroupGraphnorm}$ we have $P_t u (\cdot) = \widetilde{T_t u}(\cdot)$ is bounded uniformly in $t \ge 0$ and locally bounded in $x$. So, by Lebesgue's dominated convergence the right hand side of $\eqref{EqLaplace0}$ is continuous in $x$. 

Thus for all $x \in E_1$ 
\begin{align*}
R_\lambda u (x) = \int_0^\infty \exp(-\lambda t) P_t u (x) \, dt. 
\end{align*}

(iv): Observe that $(\lambda G_\lambda)_{\lambda > 0}$ is strongly continuous also on $(D(L_p),\Vert \cdot \Vert_{D(L_p)})$ so we get using $\eqref{EqEmbeddingCont}$ for $x \in E_1$

\begin{align}
 | \lambda R_\lambda u(x) - \widetilde u(x) | \le C_1 \Vert  \lambda G_\lambda u - u \Vert_{D(L_p)} \overset{\lambda \to \infty}{\longrightarrow} 0.
\end{align}

\end{proof}

Now we prove Theorem \ref{ClassicalStrongFeller}.
\begin{proof}[\textbf{Proof of Theorem \ref{ClassicalStrongFeller}}] \label{ProofClassicalStrongFeller}
Let $f \in \mathcal B^+_b(E)$. Set $f_n := 1_{K_n} f$, $n \in \N$, with $K_n$ as in Lemma \ref{LemmaSigmaAlgebraOpenFinite}. Then $f_n \in \mathcal L^p(E,\mu)$ and $f_n \uparrow f$. Define $u_n := R_1 f_n$, then $(1-L)u_n = f_n$ is uniformly bounded in the $L^\infty$-norm. Moreover, $(u_n)_{n \in \N}$ is also bounded in $L^\infty$-norm, since $R_1$ is sub-Markovian. So by assumption $(u_n)_{n \in \N}$ is equicontinuous. For $x \in E_1$, choose $U \subset E_1$ such that $\overline{U} \subset E_1$ is compact. Then by Arzela-Ascoli the sequence $(u_n|_{\overline{U}})_{n \in \N}$ in $C(\overline{U})$ possesses a subsequence converging uniformly to a function $v \in C(\overline{U})$. By the properties of $R_1$ and $(f_n)_{n \in \N}$ it holds that $R_1 f_n(x) \uparrow R_1 f(x)$ for every $x \in E_1$.
So $R_1 f(x) = v(x)$ for $x \in \overline{U}$. So for every $x \in E_1$ there exists a neighborhood of $x$ such that $R_1 f$ is continuous on $U$. Hence, $R_1 f \in C(E_1)$.
The claim for general $f \in \mathcal B_b(E)$ follows by linearity now.
%If $(\mathcal E,D(\mathcal E))$ is conservative, then it holds $R_1 1_E = 1$, $\mu$-a.e. Thus by the classical strong Feller property $R_1 1(x) = 1$ for all $x \in E_1$. From this we get together with the classical strong Feller property of $R_1$ that $(P_t)_{t > 0}$ is classical strong Feller as in \cite{AKR03}[Prop.~3.8].
Under the additional assumptions it follows as in \cite[Prop.~3.8]{AKR03}, that $P_t 1_E (x)= 1$ for every $t \ge 0$ and $x \in E_1$. This implies that the constructed process $\mathbf M$ is conservative.
\end{proof}

We prove an enforced version of Theorem \ref{PropKernelResolvent}(iv) which we need later on for the solution of the martingale problem.

\begin{lemma} \label{LemmaResolventLaplaceTrafo}
For all $x \in E_1$ and $f \in L^p(E,\mu) \cup \mathcal B_b(E)$ \eqref{EqResolventLaplaceTrafo} holds.
\end{lemma}
\begin{proof}
Let $x \in E_1$.
First assume that $f \in L^p(E,\mu) \cap L^\infty(E,\mu)$, $f \ge 0$ and $\Vert f \Vert_{L^\infty} \le C_5 < \infty$. 
Set $f_n := n G_n f \in D(L_p)$. Since $n G_n$ is sub-Markov we have $f_n(x) \le C_5$ $\mu$-a.e. and $f_n$ converges to $f$ in $L^p(E,\mu)$.

By \ref{PropKernelResolvent}(iii) we have
\begin{align}
R_\lambda f_n(x) = \int_0^\infty \exp(-\lambda t) P_t f_n (x) dt, \quad \text{for all} \, n \in \N. \label{EqinProofResolventLaplaceTrafo}
\end{align}

Using \eqref{EqResolventEst} and $L^p$-convergence the left hand side of $\eqref{EqinProofResolventLaplaceTrafo}$ converges to $R_\lambda f (x)$.
Furthermore, \eqref{EqSemigroupPointwise} implies for $t>0$ $\limn P_t f_n (x) = P_t f (x)$. Since $P_t$ is sub-Markovian, we get $| \exp(-\lambda t) P_t f_n (x) | \le \exp(-\lambda t) C_5$. By Lebesgue's dominated convergence the right hand side of $\eqref{EqinProofResolventLaplaceTrafo}$ converges to $\int_0^\infty \exp(-\lambda t) P_t f (x) dt$.
Thus, $\eqref{EqResolventLaplaceTrafo}$ holds for those $f$.\\
Now let $f \in L^p(E,\mu)$ with $f \ge 0$ set $f_n := f \wedge n$. Then $f_n \uparrow f$ and $\eqref{EqinProofResolventLaplaceTrafo}$ holds for all $f_n$. By monotone convergence on both sides we get the identity for $f$.

Now for $f \in L^p(E,\mu)$ observe that $P_t f (x) = P_t f^+ (x) - P_t f^- (x)$ and therefore $|P_t f(x)| \le P_t f^+ (x) + P_t f^- (x)$. 
Since by the proven statement $\int_0^\infty \exp(-\lambda t) P_t f^{+/-} (x) dt < \infty$ we get that the integral  $\int_0^\infty \exp(-\lambda t) P_t f (x) dt$ exists. Then $\eqref{EqinProofResolventLaplaceTrafo}$ follows by linearity of $R_\lambda$ and $P_t$.
Thus, the class of all functions satisfying $\eqref{EqinProofResolventLaplaceTrafo}$ fulfills \ref{TheoFMCT}(i), \ref{TheoFMCT}(iii) and \ref{TheoFMCT}(ii'). So the statement for $f \in \mathcal B_b(E)$ follows using Corollary \ref{CoroTheoFMCT}.
\end{proof}

Based on this lemma we prove a pointwise equation relating the semigroup of kernels and the generator. This formula is essential for the solution of the martingale problem.

\begin{lemma} \label{LemmaSemigroupDiff}
For $x \in E_1$, $u \in D(L_p)$ it holds for all $t > 0$
\begin{align}
P_t u(x) - \widetilde u(x) = \int_0^t P_s L_p u(x) ds = \int_0^t \widetilde {L_p P_s u}(x) ds \label{EqLemmaSemigroupDiff}
\end{align}
and the integral is well-defined.
Here $\widetilde u(x)$ ($\widetilde {L_p P_s u} (x)$) denotes the value of the respective continuous version at $x$.

\end{lemma}
\begin{proof}
Let $x \in E_1$, $t > 0$ be fixed.
First note that for $f \in L^p(E,\mu)$ the map $[0,t] \ni s \mapsto P_s f(x) \in \R$ is integrable. Indeed, consider $f \ge 0$ first. Then
\begin{multline*}
\int_0^t P_s f(x) ds = \int_0^t \exp(s) \exp(-s) P_s f(x) ds \\
\le \exp(t) \int_0^t \exp(-s) P_s f(x) ds \le \exp(t) R_1 f(x) < \infty.
\end{multline*}
Note that the integral value is independent of the $\mu$-version of $f$.
Since $|P_s f| \le P_s |f|$ the statement follows for general $f \in L^p(E,\mu)$.

Let $u \in D(L_p)$, set $f := (1 -L_p) u$, then $u = G_1 f$. Since $T_t G_1 f = G_1 T_t f$, we have by continuity $P_t R_1 f(x) = R_1 P_t f(x)$. Using Lemma \ref{LemmaResolventLaplaceTrafo} and the semigroup property we get
\begin{multline*}
\exp(-t) P_t R_1 f(x) =  R_1 \exp(-t) P_t f \\
= \int_0^\infty \exp(-(t+s)) \, P_{t+s} f(x) \, ds
= \int_t^\infty \exp(-s) \, P_s f(x) \, ds.
\end{multline*}
Thus,
\begin{multline*}
\exp(-t) P_t R_1 f(x) - R_1 f(x) = - \int_0^t \exp(-s) P_s f(x) ds
= - \int_0^t \exp(-s) P_s (1 - L) u(x) ds. 
\end{multline*}

By construction $R_1 f$ is the unique continuous version of $G_1 f$, thus
\begin{align}
\exp(-t) P_t u(x) - \tilde u(x) = - \int_0^t \exp(-s) P_s (1 - L) u(x) ds. \label{EqExpSemigroupDiff}
\end{align}

So for every $x \in E_1$ the mapping $[0,t] \ni s \mapsto \exp(-s) P_s \tilde u(x)$ is absolutely continuous with integrable weak derivative $\exp(-s) P_s (L - 1) u(x)$, $s>0$. So by the product rule also the mapping $s \mapsto P_s \tilde u(x)$ is absolutely continuous with weak derivative $P_s L u(x)$, $s > 0$. This proves the first equality of \eqref{EqLemmaSemigroupDiff}. For the second equality note that $L_p T_s u = T_s L_p u$ for $s > 0$ and since $T_s u \in D(L^2_p)$ we have by continuity that $\widetilde {L_p P_s u}(x) = P_s L_p u(x)$ for every $x \in E_1$. 
\end{proof}

\section{Construction of the $\mathcal L^p$-strong Feller process} \label{SectionProcess}
In this section we apply the results from the previous section to construct the corresponding diffusion process. Throughout this section the same assumptions as in Section \ref{SectionConstructionKernel} are assumed. The construction is based on the techniques developed in \cite{AKR03} and \cite{Doh05}.

The process is first constructed with time parameter in the dyadic numbers. Using the diffusion process which is associated to a regular strongly local Dirichlet form, we derive properties of the constructed process. The essential step is to get from quasi-everywhere statements to pointwise statements. For this we need a probability measure on $E$ which has the same nullsets as $\mu$.

\begin{lemma}
Let $(S,\mathcal B,\mu)$ be $\sigma$-finite with $\mu(S)>0$. Then there exists a probability measure $\nu$, i.e. $\nu(E) = 1$, absolutely continuous to $\mu$ such that $\nu(A) = 0$ implies $\mu(A)=0$.
\end{lemma}
\begin{proof}
Choose an increasing sequence of measurable sets $F_k \in \mathcal B$ with $\mu(F_k) < \infty$, $k \in \N$, and  $\bigcup_{k \in \N} F_k = S$. Define $G_1 := F_1$ and inductively $G_{k+1}:= F_{k+1} \setminus F_k$, $k \in \N$. 
Set $\alpha_k := 2^{-k} \frac{1}{\mu(G_{k})} < \infty$ if $\mu(G_{k}) > 0$ and $\alpha_k = 0$ otherwise.
Define $\alpha = \sum_{k=1}^\infty \alpha_k 1_{G_k}$. Then $\alpha > 0, \mu$-a.e.,~and $0 < \int_{S} \alpha d \mu < \infty$.
Now normalize $\alpha$ and set $\nu := \alpha \, \mu $. If $\nu(A)=0$ for some $A \in \mathcal B$, then for $G_k$ with $\mu(G_k) \ne 0$ it holds $\int_{G_k \cap A} \alpha \, \mu = 0$. Since $\alpha>0, \mu$-a.e.,~on $G_k$ it follows $\mu(G_k \cap A)=0$. So by construction of the $G_k$ it follows $\mu(A)=0$.
\end{proof}

Now we fix such a measure $\nu$ for $(E,\mathcal B(E),\mu)$. Thus, we get a probability distribution of starting points that has the same nullsets as $\mu$.

We extend $E$ by the cemetery point $\Delta$ and endow $E^\Delta := E \cup \{ \Delta \}$ with the topology of the Alexandrov compactification, i.e., $E^\Delta$ is the one-point compactification of $E$. Note that there exists a complete metric on $E^\Delta$ inducing the Alexandrov topology, see e.g. \cite[Cor.~3.45]{CB06}.
Next we extend the semigroup of kernels $(P_t)_{t \ge 0}$ from $E_1 \times \mathcal B (E)$ to $E^\Delta \times \mathcal B(E^\Delta)$. Here $\mathcal B(E^\Delta) = \mathcal B(E) \cup \{ \mathcal B(E) \cup {\Delta} \}$.

\begin{definition} \label{DefinitionPDeltaKernel}
Let $(P_t)_{t \ge 0}$ be the kernels on $E_1 \times \mathcal B (E)$ from Theorem \ref{TheoremKernelSemigroup}.
Define $P^\Delta_t : E^\Delta \times \mathcal B(E^\Delta)$ by:
If $x \in E^\Delta \setminus E_1$, $t \ge 0$ then $P^\Delta_t(x,\cdot):= \varepsilon_x$. Here $\varepsilon_x$ is the point measure in $x$.
For $x \in E_1$, $t > 0$, $A \in \mathcal B(E^\Delta)$ define $P^\Delta_t(x,A) = P_t(x,A \setminus \{ \Delta \}) + (1-P_t(x,E)) \varepsilon_{\Delta}(A)$.  
For $t=0$ define $P^\Delta_t(x,\cdot) := \varepsilon_x$.
\end{definition}
A straightforward calculation gives that $(P^\Delta_t)_{t \ge 0}$ is a semigroup of kernels.
The measure $P^\Delta_t(x,\cdot)$ consists of the part $P_t(x,\cdot)$ which is absolutely continuous w.r.t $\mu$ and the singular part $\varepsilon_{\Delta}$.
We extend each $f \in \mathcal B(E)$ to $\mathcal B(E^\Delta)$ by $f(\Delta) := 0$. For $n \in \N$, we set
\begin{align*}
S_n := \{ k 2^{-n} \, | \, k \in \N \cup \{0\} \} \ \text{and} \ S:= \bigcup_{m \in \N} S_m.
\end{align*}
$(P^\Delta_t(x,dy))_{t \ge 0}$ is a semigroup of probability kernels on the polish space $E^\Delta$, so by Kolmogorov's standard construction scheme, see e.g. \cite[Ch.~I, Theo.~2.11]{BG68}, there exists a family of probability measures $\mathbf P_x$, $x \in E^\Delta$, on $\mathbf \Om := (E^\Delta)^S$, equipped with the product $\sigma$-field $\mathcal F^0$, such that
\begin{align}
\mathbf M^0 := \left(\mathbf \Om^\Delta,\mathcal F^0, (\mathcal F^0_s)_{s \in S},(\mathbf X_s^0)_{s \in S},(\mathbf P_x)_{x \in E^\Delta} \right)
\end{align}
is a normal Markov process having transition kernels $(P^\Delta_t(x,\cdot))_{t \ge 0}$.
Here $\mathbf X^0_s : \mathbf \Om \to E^\Delta$ are the coordinate maps and $\mathcal F^0_s := \sigma(\mathbf X^0_t \, | \, t \le s)$.

This process is defined only for dyadic time parameters at first. Next we show that this process can be uniquely extended to a process with time parameter $t \in [0,\infty)$ having continuous paths on $(0,\infty)$. Under Condition \ref{CondRegularityDomain}(ii) this process has even right continuous paths at $t=0$. To study the properties of $\mathbf M^0$ we use the result of \cite[Theo.~4.5.3]{FOT94}, see also \cite[Ch. V, Theo.~1.11]{MR92}. Here and in the sequel of this article $\mathcal E$-quasi statements hold outside some set $S \in \mathcal B$ with cap$_{\mathcal E}(S)$=0, i.e.,~capacity zero. There exists a diffusion (i.e., a strong Markov process having continuous sample paths) 
\begin{align*}
\mathbf {\mhat{M}} = \left(\mathbf{\mhat{\Om}},\mathcal{ \mhat{F}},(\mathcal {\mhat{F}}_t)_{t \ge 0},(\mathbf{\mhat{X}}_t)_{t \ge 0},(\mhat{\mathbf P}_x)_{x \in E^\Delta} \right),
\end{align*}
which is properly associated with the regular Dirichlet form $(\mathcal E,D(\mathcal E))$, i.e.,  $\mhat{P}_t f$ is a $\mathcal E$-quasi continuous version of $T_t f$ for $f \in L^2(E,\mu) \cap \mathcal B_b(E)$, where $\mhat{P}_t f (x) = \mathbb E_x[f(\mathbf{\mhat X}_t)]$, see \cite[Ch. III, Def.~2.5]{MR92}. 
For a Borel set $K$ define $\sigma(K) = \text{inf} \{ t > 0 \ | \ \mathbf{\mhat{X}}_t \in K \}$, the first hitting time. Moreover, define the lifetime $\lifetime := \text{inf} \{ t > 0 \ | \ \mathbf{\mhat X}_t=\Delta\}$.
The strong local property of $(\mathcal E,D(\mathcal E))$ implies that $\mathbf{\mhat{M}}$ enters the cemetery only continuously, see \cite[Theo.~4.5.3]{FOT94}. So $(\mathbf{\mhat{X}}_t)_{t \ge 0}$ is even continuous for all $t \in [0,\infty)$ and not only for $t \in [0,\lifetime)$. Thus $\mathbf{\mhat \Om}$ may be chosen as $C([0,\infty), E^\Delta)$.

\begin{lemma} \label{CapZeroCompSet}
Assume that for $E_1 \subset E$ it holds $\text{cap}_{\mathcal E}(E \setminus E_1)$=0. Then there exists a sequence of compact sets $(K_n)_{n \in \N}$ with $K_n \subset E_1$ such that 
\begin{align}
\mathbf{\mhat P}_x \left(\limn \sigma_{K^c_n} \ge \lifetime \right) = 1 \ \ \text{for quasi-every} \ \ x \in E, \label{EqPKn}
\end{align}
in particular, $\mathbf{\mhat P}_\nu \left(\limn \sigma_{K^c_n} \ge \lifetime \right) = 1$, where $\mathbf{\mhat P}_\nu := \int_E \mathbf{\mhat P}_x \nu(dx)$.
\end{lemma}
We use the notation $K^c = E \setminus K$ for the complement of $K$ in $E$.
\begin{proof}
From the definition of the capacity and Lemma \ref{LemmaSigmaAlgebraOpenFinite} there exists a sequence of compacts sets $(K_n)_{n \in \N} \subset E_1$ with the property that $\limn \text{cap}_{\mathcal E}(K \setminus K_n) = 0$. Then \eqref{EqPKn} follows from \cite[Ch.~IV, Lem.~4.5]{MR92} or \cite[Lem.~5.1.6]{FOT94}.
\end{proof}
%\begin{proof}
%Since $(E,d)$ is locally compact and separable, there exist a sequence $(C_n)_{n \in \N}$ of compact set with $C_n \subset \interior{C}_{n+1}$ and $E= \cup_{n \in \N} C_n$. Moreover, there exist a sequence $(U_n)_{n \in \N}$ of open sets with $E_0 \subset \cap_{n \in \N} U_n$ and $\limn cap_{\mathcal E}(U_n)=0$. Define $K_n := C_n \cap U^c_n$. Since $U^c_n$ is closed and $C_n$ is compact, $K_n$ is compact. By construction $K_n \subset E^c_0$. Now let $K \subset E$ be an arbitrary compact set. Since $K$ is compact and $(\interior{C}_n)_{n \in \N}$ is an open cover, there exists finitely many $C_1,...,C_m$ covering $K$. So for $l > m$ we have that $K \cap C^c_l = \emptyset$. Thus $K \setminus K_l = K \cap (U_l \cup C^c_l) = K \cap U_l$.
%Thus 
%\begin{align*}
%cap_{\mathcal E}(K \setminus K_l) \le cap_{\mathcal E}(U_l) \underset{l \to \infty}{\longrightarrow} 0.
%\end{align*}
%By \cite{MR92}[Lemma IV.4.5] or \cite{FOT94}[Lemma 5.1.6] it follows that $\mathbb P_x(\limn \sigma_{K^c_n} \ge \lifetime) = 1$ for quasi-every $x \in E$. Now note that this equality holds also for $\mu$-a.e. $x \in E$. Thus the second conclusion follows, since $\nu$ is absolutely continuous to $\mu$.
%\end{proof}

In the sequel we fix one sequence of compact sets  $(K_n)_{n \in \N}$ with the properties as in Lemma \ref{CapZeroCompSet}.
\begin{remark} \label{RemCapZero}
Every sequence of closed sets $(F_n)_{n \in \N} \subset E$ with the property that \newline $\limn \text{cap}_{\mathcal E}(K \setminus F_n) = 0$ for every compact set $K \subset E$ has property $\eqref{EqPKn}$. Nevertheless, in concrete applications it might be more suitable to construct an explicit sequence of compact sets with this property. All further constructions work for such a sequence of sets, not only for the sets provided by Lemma \ref{CapZeroCompSet}.
\end{remark}

Define
\begin{align}
\mathbf{\mhat \Om}_0 := \left\{ \omega \in \mathbf{\mhat \Om} \ | \ \omega(0) \in E_1 \cup \{ \Delta \}, \underset{n \to \infty}{\text{lim}} \sigma_{K^c_n} \ge \lifetime, \omega (t) = \Delta, t \ge \lifetime \right \},
\end{align}
then it holds
\begin{align*}
\mathbf{\mhat P}_x(\mathbf{\mhat \Om}_0) = 1 \ \ \text{for quasi-every} \ \ x \in E_1 \cup \{ \Delta \}.
\end{align*}
Moreover, for
\begin{align}
\mathbf{\mhat P}_{\nu} := \int_{E} \mathbf{\mhat P}_x \nu(dx) \label{EqTildePnu}
\end{align}
it holds
\begin{align*}
\mathbf{\mhat P}_{\nu}(\mathbf{\mhat \Om}_0) = 1,
\end{align*}
i.e., with $\mathbf{\mhat P}_{\nu}$-probability $1$ we observe continuous paths starting from $E_1 \cup \{ \Delta \}$ which stay in $E_1 \cup \{ \Delta \}$ and reach $\Delta$ only continuously.
Define the map $G: \mathbf{\mhat \Om} \to \mathbf \Om$ by
\begin{align*}
\mathbf{\mhat \Om} \ni \omega = (\omega(t))_{t \in [0,\infty)} \mapsto G(\omega) := (\omega(s))_{s \in S}.
\end{align*}
Since every continuous function is uniquely determined by its values on a dense set, $G$ is a one-to-one map. Moreover, for $\mathbf{\mhat{\mathcal F}}^0 := \sigma ( \mathbf{\mhat X}_s \, | \, s \in S )$ it holds $\mathbf{\mhat \Om}_0 \in \mathbf{\mhat{\mathcal F}}^0$ and $G$ is $\mathbf{\mhat{\mathcal F}}^0 / \mathbf{\mathcal F}^0$-measurable, see \cite{Doh05}.

With the following lemma we can connect the measures $\mathbf P_\nu$ and $\mathbf {\mhat P}_\nu$.
\begin{lemma} \label{LemmaImageMeasurePhead}
Define $\mathbf{\mhat P}'_\nu$ as the image measure of the restriction of $\mathbf{\mhat P}_{\nu}$ to $\mathbf{\mhat{\mathcal F^0}}$ under $G$, i.e.
\begin{align*}
\mathbf{\hat P}'_\nu := \mathbf{\mhat P}_{\nu}|_{\mathbf{\mhat{\mathcal F}}^0} \circ G^{-1}.
\end{align*}
Then
\begin{align*}
\mathbf{\mhat P}'_\nu = \mathbf P_\nu,
\end{align*}
where $\mathbf P_\nu$ is defined as in \eqref{EqTildePnu} with $\mathbf{\mhat P}_x$ replaced by $\mathbf P_x$.
\end{lemma}
The proof works as in \cite[Lem.~4.2]{AKR03} by considering sets $K_0,...,K_k \subset \mathcal B(E)$ of finite measure and time points $0 < t_1 < ... < t_k$, $k \in \N$, and using that both ${\hat P}_{t_k} f$ and $P_{t_k} f$ are $\mu$-version of $T_{t_k} f$ for a bounded function $f$ with support having finite $\mu$-measure.
Since $\mathcal B(E)$ is generated by those sets and both transition kernels are extended in the same way to $E^\Delta$ the claim follows.

As in \cite[Lem.~4.3]{AKR03} we get using Lemma \ref{LemmaImageMeasurePhead} the following Lemma.
\begin{lemma} \label{LemmaPmu1}
$G(\mhat{\mathbf \Om}_0) \in \mathbf{\mhat{\mathcal F}}^0$ and $\mathbf P_y(G(\mathbf{\mhat \Om}_0)) = 1$ for $\mu$-a.e. $y \in E_1$.
\end{lemma}

Consider the time shift operator $\theta_s : \mathbf \Om \to \mathbf \Om$, $\theta_s(\omega) = \omega(\cdot + s)$, $s \ge 0$, and define
\begin{align*}
\mathbf \Om_1 := \bigcap_{s > 0,s \in S} \theta^{-1}_s(G(\mathbf{\mhat \Om}_0)),
\end{align*}
i.e., all paths with time parameter in $S$ which come from a path on $[0,\infty)$ which is continuous in $(0,\infty)$ and does not hit $E \setminus E_1$.
Then using Lemma \ref{LemmaPmu1} we get
\begin{lemma} \label{LemmaOmega1P1}
\begin{align*}
\mathbf P_x(\mathbf \Om_1) = 1 \quad \text{for} \ x \, \in \, E_1 \cup \{ \Delta \}.
\end{align*}
\end{lemma}
\begin{proof}
Follows using the Markov property of $\mathbf M^0$ and the absolute continuity of $(P^\Delta_t)_{t > 0}$ on $E$ and the fact that $P_s^\Delta(\Delta,G(\mathbf{\mhat \Om}_0))=1$, see e.g. \cite[Lem.~4.4]{AKR03}.
\end{proof}
%\begin{proof}
%We prove that for every $s > 0$, $s \in S$, $y \in E_1$ it holds 
%\begin{align*}
%\mathbf P_y (\theta_s^{-1}(G(\tilde{\mathbf \Om}_0))) = 1.
%\end{align*}
%
%Using the Markov property of $\mathbf M^0$ in $*$ we get
%\begin{multline*}
%\mathbf P_x(\theta^{-1}_s(G(\tilde{\mathbf \Om}_0))) = \mathbf P_x( (X_{t+s,\, t \ge 0 ,\, t \in S}) \in G(\tilde{\mathbf \Om}_0)) = \mathbb E_x[ \mathbb E_x[1_{X_{t+s,\, t \ge 0 ,\, t \in S} \in G(\tilde{\mathbf \Om}_0)} | \mathcal F^0_s]] \\
%\overset{*}{=} \mathbb E_x[\mathbf P_{X_s}(G(\tilde{\mathbf \Om}_0))] 
%= \int_{E} p_s(x,y) \mathbf P_{y}(G(\tilde{\mathbf \Om}_0)) \mu(dy) + (1 - P_s(x,E)) \mathbf P_\Delta(G(\tilde{\mathbf \Om}_0)).
%\end{multline*}
%By definition of $\mathbf P_\Delta$ and $\tilde {\mathbf \Om}_0$ we have $\mathbf P_\Delta(G(\tilde {\mathbf \Om}_0)) = 1$ and by \ref{LemmaPmu1} we have for $\mu$-a.e. $y \in E$ $\mathbf P_y(G(\tilde {\mathbf \Om}_0))=1$.
%Thus
%\begin{align*}
%\mathbf P_x(\theta^{-1}_s(G(\tilde{\mathbf \Om}_0))) = 1. 
%\end{align*}
%
%But then we also have
%\begin{align*}
%\mathbf P_x \left(\bigcap_{s > 0,s \in S} \theta^{-1}_s(G(\tilde{\mathbf \Om_0})) \right) = 1.
%\end{align*}
%\end{proof}

To get right continuity of the paths at $t=0$ we use the point separating Condition \ref{CondRegularityDomain}(ii).

\begin{lemma} \label{LemmaContTZero}
For $x \in E_1$ it holds
\begin{align}
\underset{s \downarrow 0, \, s \in S}{\lim} \mathbf X^0_s = x \quad \mathbf P_x\text{-a.s.} \label{EqContXatZero}
\end{align}
\end{lemma}
\begin{proof}
Let $x \in E_1$.
From Lemma \ref{LemmaResolventLaplaceTrafo} we get that for $f \ge 0$, $(\exp(-s) R_1 f(\mathbf X^0_s))_{s \in S}$ is an $\mathcal F^0_s$-supermartingale.

Now choose the point separating functions $(u_n)_{n \in \N} \subset D(L_p)$ from Condition \ref{CondRegularityDomain}.
Then as in \cite[Lem.~4.6]{AKR03} we obtain using supermartingale convergence and Theorem \ref{TheoremKernelSemigroup}(iii), that $\mathbf P_x$-a.s.
\begin{align}
\underset{s \downarrow 0,\, s \in S}{\lim} u_n(\mathbf X^0_s) = u_n(x), \quad \, \text{for all} \  n \in \N. \label{EqLimUk}
\end{align}
Using the point separating property of $(u_n)_{n \in \N}$ this implies
\begin{align*}
\underset{s \downarrow 0,\, s \in S}{\lim} \mathbf X^0_s = x.
\end{align*}
\end{proof}

After these preparations we can construct the $\mathcal L^p$-strong Feller process.

\begin{proof}[\textbf{Construction of the process of Theorem \ref{TheoProcess}}] \label{proofTheoDiffProcess}
Choose a fixed $x_0 \in E_1$ and define for $\omega \in \mathbf \Om$, $t \ge 0$
\begin{align*}
\mathbf X_t(\omega) :=
\begin{cases}
\lim_{s \downarrow t,s \in S} \mathbf X_s(\omega) & \ \text{if} \ \omega \in \mathbf \Om_1 \, \text{and the limit exists}, \\
x_0 & \text{else}.
\end{cases}
\end{align*}
Furthermore, let $(\mathbf{\mathcal F}_t)_{t \in [0,\infty]}$ be the corresponding natural filtration as in \cite[Ch.~IV, Def.~1.8,(1.6) and (1.7)]{MR92}.\\
First let $t>0$ arbitrary.
Then by Lemma \ref{LemmaOmega1P1} we have for $x \in E_1$ that $\mathbf X_t$ is $\mathbf P_x$-a.s. the limit of $(\mathbf X_{s_n})_{n \in \N}$ for every dyadic sequence $(s_n)_{n \in \N}$ with $s_n \downarrow t$.
Moreover the \textit{law} $\mathcal L^x(X_{s_n})$ of $X_{s_n}$ is given by $P^\Delta_{s_n}(x,\cdot)$. 
Now let $u \in C_b(E) \cap \mathcal L^p(\Om,\mu)$. Then using Lebesgue's dominated convergence in the first equality and Theorem \ref{TheoremKernelSemigroup}$(iii)$ in the third one we get
\begin{align*}
\mathbb E_x[u(\mathbf X_t)] = \limn \mathbb E_x[u(\mathbf X_{s_n})] = \limn P_{s_n} u (x) = P_t u(x).
\end{align*}
Corollary \ref{CoroTheoFMCT} implies that 
\begin{align}
\mathbb E_x[u(\mathbf X_t)] = P_t u(x) \label{EqProofLawXt}
\end{align}
holds for $u \in \mathcal B_b(E)$. Furthermore, note that $\mathbf P_x(\mathbf X_t = \Delta) = 1 - \mathbf P_x(\mathbf X_t \in E) = 1- P_t 1_E (x)$.
Thus $\mathcal L_x(\mathbf X_t) = P^\Delta_t(x,\cdot)$.

For $t=0$ we have by Lemma \ref{LemmaContTZero} that $(\mathbf X_{s_n})_{n \in \N}$ converges to $x$ for every dyadic sequence $(s_n)_{n \in \N}$ with $s_n \downarrow 0$. So $\mathcal L((\mathbf X_t)_{t \ge 0}) = (P^\Delta_t)_{t \ge 0}$.\\ Together with Lemma \ref{LemmaOmega1P1}, $(\mathbf X_t)_{t \ge 0}$ has $\mathbf P_x$-a.s.~continuous sample paths with $\mathbf X_0 = x$.\\
%\end{proof}
%\begin{remark}
%From the previous theorem we conclude using Lebesgue's dominated convergence on $(\mathbf \Om,\mathbf P_x)$, $x \in E_1$, that for all functions $f \in C_b(E^\Delta)$ it holds
%\begin{align*}
%P^\Delta_s f (x) \overset{s \to t} \longrightarrow P^\Delta_t f (x) \quad \text{for} \ x \in E_1. 
%\end{align*}
%\end{remark}

%Using Corollary \ref{CoroTheoFMCT} and the $\mathcal L^p$-strong Feller property of the resolvent the strong Markov %property follows similarly as in \cite[Ch.~I, Theo.~8.11]{BG68}. 

The strong Markov property follows as in \cite[Ch.~I,Theo.~8.11]{BG68}. There as set $\mathbf L \subset C_b(E_1)$ we may choose $C_b(E_1) \cap \mathcal L^p(E,\mu)$ in our case, using Corollary \ref{CoroTheoFMCT} and the $\mathcal L^p$-strong Feller property of the resolvent.

\end{proof}

\begin{lemma} \label{LResolventProc}
For $x \in E_1$, $\lambda > 0$, $u \in \mathcal L^p(E,\mu) \cup \mathcal B_b(E)$ it holds
\begin{align*}
\mathbb E_{x} \left[\int_0^\infty \exp(-\lambda s) u(\mathbf X_s) ds \right] = R_\lambda u(x). 
\end{align*}
\end{lemma}
\begin{proof}
%Let $u \in \mathcal L^p(E,\mu)$, $f \ge 0$.
%From Lemma \ref{LemmaResolventLaplaceTrafo} and the fact that $\mathcal L_x (X_t) = P_t(x,\cdot)$ we get
%\begin{align*}
% \int_0^\infty \int_{\mathbf \Om} \exp(-\lambda s) u(\mathbf X_s(\omega)) d \mathbf P_x(\omega) ds = \int_0^\infty \exp(-\lambda s) P_s u(x) ds = R_\lambda u(x) < \infty.
%\end{align*}
%So by Fubini-Tonelli it follows
%\begin{multline*}
% \mathbb E_x \left [ \int_0^\infty \exp(-\lambda s) u(X_s) ds \right ] = \int_{\mathbf \Omega} \int_0^\infty \exp(-\lambda s) u(\mathbf X_s(\omega)) ds d \mathbf P_x(\omega) \\
%= \int_0^\infty \int_{\mathbf \Omega} \exp(-\lambda s) u(\mathbf X_s(\omega)) d \mathbf P_x(\omega) ds = R_\lambda u(x). 
%\end{multline*}
%
%The statement for general $u \in \mathcal L^p(E,\mu)$ follows by applying this to $|u|$, then to $u^+$ and $u^-$ and then using linearity.
%For $u \in \mathcal B_b(E)$ the statement follows using Corollary \ref{CoroTheoFMCT}.
Follows by Lemma \ref{LemmaResolventLaplaceTrafo} and $\mathcal L_x(\mathbf X_t) = P^\Delta_t(x,\cdot)$, $t \ge 0$, $x \in E_1$, using Fubini and Corollary \ref{CoroTheoFMCT}.
\end{proof}

Next we state some properties of the process $\mathbf M$ which can be transferred from the process $\mathbf{\mhat M}$ to pointwise statements on $E_1$.

\begin{theorem} \label{TheoPropHit}
Let $\mathbf M=(\mathbf \Om,\mathbf{\mathcal F},(\mathbf{\mathcal F}_t)_{t \ge 0},(\mathbf X_t)_{t \ge 0},(\mathbf P_x)_{x \in E_1 \cup \{ \Delta \} })$ be the Markov process from Theorem \ref{TheoProcess}. Then the following properties hold:

\begin{enumerate}[label=(\roman{*}),ref=(\roman{*}]
 \item If $U \subset E_1$ has $\mu$-measure zero, then $dx(\{ s \in \R_0^+ \, | \, \mathbf X_s \in U \}) = 0$, $dx$ the Lebesgue measure on $\R_0^+$, $\mathbf P_x$-a.s. for $x \in E_1$.
 \item If $U \subset \mathcal B(E)$ has the property  $cap_{\mathcal E}(U) = 0$, then $\mathbf P_x(\sigma_U < \infty) = 0$ for $x \in E_1$.
% \item If $R_1 1 (x)= 1$ for $x \in E_1$, then $\mathbf M$ is conservative.
% \item If $(\mathcal E,D(\mathcal E))$ is conservative and $(R_\lambda)_{\lambda > 0}$ is strong Feller in the classical sense, then $\mathbf M$ is conservative. 
\end{enumerate}
\end{theorem}

\begin{proof}
$(i)$: If $\mu(U) = 0$ then $G_1 1_U(x) = 0$ for $\mu$-a.e.~$x \in E$. By the $\mathcal L^p$-strong Feller property it holds $R_1 1_U (x) = 0$ for $x \in E_1$. Thus by Lemma \ref{LResolventProc}
\begin{align*}
 \mathbf E_x \left[\int_0^\infty \exp(-t) 1_U(\mathbf X_t) dt \right] = R_1 1_U (x) = 0.
\end{align*}
Hence $\mathbf P_x$-a.s. it holds that $\exp(-t) 1_U(\mathbf X_t) = 0$ for $dx$-a.e.~$t$, hence $1_U(\mathbf X_t) = 0$ for $dx$-a.e.~$t$. \newline
$(ii)$: 
%Assume that $U \in \mathcal B(E)$ and $cap_{\mathcal E}(U)=0$. Then it holds for $q.e.$ $x \in E$ \newline
%that $\mathbf P_x(\sigma_U = \infty)$, see e.g. \cite{FOT94}[Theo.~4.2.1] \footnote{The additional assumption there that every compact set has finite capacity is only used in the converse direction of the implication.}. As in the proof of \ref{LemmaOmega1P1} we get that for $x \in E_1$ it holds $\mathbf P_x(\exists \, s \ge t \, \text{such that} \, X_s \in U) = 0$ for every $t > 0$. Thus $\mathbf P_x(\exists \, s > 0 \, \text{such that} \, X_s \in U)=0$. Hence $\sigma_U = \infty$, $\mathbf P_x$-a.s. for $x \in E_1$. \newline
Follows as in the preparatory Lemmas for the proof of Theorem $\ref{TheoProcess}$.

%\noindent $(iii)$ 
%We have by Lemma \ref{LResolventProc}
%\begin{align*}
% \mathbb E_x \left [ \int_0^\infty \exp(- s) 1_E (\mathbf X_s) ds \right ] = R_1 1(x)= 1.
%\end{align*}
%So $\mathbf P_x$-a.s. it holds $ds$-a.s. $1_E(X_s(\omega)) = 1$, i.e. $X_s(\omega) \in E$. 
%Now assume there exists $s \in [0,\infty)$ such that $X_s(\omega) = \Delta$. Then $X_t(\omega) = \Delta$ for $t > s$. This is a contradiction. \newline
%$(iv)$ 
%From conservativity of $(\mathcal E,D(\mathcal E))$ we get $R_1 1_E = 1$, $\mu$-a.e., hence by $\mathcal B_b$-strong Feller property $R_1 1_E(x) = 1$ for every $x \in E_1$. Now apply $(iii)$.
%We have for $t >0$, $T_t 1 (x) = 1$ for $\mu$-a.e. $x \in E$. Hence also $P_t 1(x)=1$ for $\mu$-a.e. $x \in E_1$. 
%Thus we get
%\begin{align*}
% R_1 1(x) = \int_0^\infty \exp(-s) p_s 1(x) ds = 1,
%\end{align*}
%for $\mu$-a.e. $x \in E_1$. This holds, since the integral can be approximated by a limit of sums evaluating $p_s 1(x)$ only at a countable number of timepoints. Together with the classical strong Feller property we get $R_1 1(x) = 1$ for every $x \in E_1$. Thus by $(iii)$, $\mathbf M$ is conservative.
\end{proof}

%\begin{remark}
%If one drops in \ref{TheoPropHit}.$(ii)$ the assumption $U \in \mathcal B(E)$ then the conclusion is weakened in the following form: There exists $U_n \subset E$, $n \in \N$, open with $U \subset \cap_{n \in \N} U_n$ and $\sigma_{\cap_{n \in \N}} U_n = \infty$, $\mathbf P_x$-a.s. for $x \in E_1$. The existence of such sets readily follows by the definition of the capacity. 
%\end{remark}

%\section{Solution of the martingale problem} \label{SectionMartingale}

We prove that the process solves the martingale problem for functions $(L_p,D(L_p))$.
\begin{theorem} \label{TheoremMartingaleSolution}
The diffusion process $\mathbf M=(\mathbf \Om,\mathcal F,(\mathcal F_t)_{t \ge 0},(\mathbf X_t)_{t \ge 0},(\mathbf P_x)_{x \in E_1 \cup \{ \Delta \}})$ of Theorem \ref{TheoProcess} solves the martingale problem for $(L_p,D(L_p))$, i.e.,
\begin{align*}
M_t^{[u]}:=\widetilde u({\mathbf X}_t) - \widetilde u(x) - \int_0^t{L_p u({\mathbf X}_s)~ds},~t\geq0,
\end{align*}
is an $(\mathcal{F}_t)$-martingale under $\mathbf{P}_x$ for all $u \in D(L_p)$, $x \in E_1$.
%, i.e.
%\begin{align*}
%M^{[u]}_t := u(\mathbf X_t) - u(\mathbf X_0) - \int_0^t L_p u(\mathbf X_s) ds,\ t \ge 0,
%\end{align*}
%is an $(\mathcal F_t)_{t \ge 0}$-martingale starting at zero. Here $u$ is choosen to be the unique continuous version.
\end{theorem}
\begin{proof}
Let $x \in E_1$, $u \in D(L_p)$. Then by Lemma \ref{LemmaSemigroupDiff} it holds
\begin{align*}
P_t u(x) - \widetilde u(x) = \int_0^t P_s L_p u (x)ds
\end{align*}
and the integral is well-defined.
Since both $u$ and $L_p u$ are extended by $0$ to $E^\Delta$, we have that this equality also holds with $P_t$ and $P_s$ replaced by $P^\Delta_t$ and $P^\Delta_s$ respectively. 
Since $\int_0^t \mathbb E_x[(L_p u)^{+/-}(\mathbf X_s)] ds = \int_0^t P_s (L_p u)^{+/-} (x) ds < \infty $ the integral $\int_0^t L_p u(\mathbf X_s) ds$ exists $\mathbf P_x$-a.s. 
Note that by Theorem \ref{TheoPropHit}(i) the integral is independent of the $\mu$-version of $L_p u$.
Now the statement follows using the Markov property of $(\mathbf X_t)_{t \ge 0}$ and that $(P^\Delta_t)_{t \ge 0}$ is the transition semigroup of $(\mathbf X_t)_{t \ge 0}$.
\end{proof}

%\section{Some remarks on the continuity}
%We give some supplementary remarks concerning the different notions of continuity.
%It was assumed that for $u \in D(L_p)$ it holds that $u \in C(E_1)$. This can be checked in the following way:
%\begin{lemma}
% If for $x \in E_1$ there exists a neighborhood $U$ of $x$ with respect to the trace topology on $E_1$ or with respect to the topology on $E$ such that $u_{|U}$ is continuous, then $u \in C(E_1)$.
%\end{lemma}
%\begin{proof}
%Since $(E_1,d|_{E_1})$ is a metric space, it suffices to prove sequentially continuity. Let $(x_n)_{n \in \N} \subset E_1$, $U$ a neighborhood of $x$ in $E_1$. Then there exists $N \in \N$ such that $x_m \in U$ for $m > \N$. So
%\begin{align*}
% \limn u(x_n) = \limn u_{|U}(x_n) = \limn u(x).
%\end{align*}
%If $U$ is a neighborhood of $x$ in $E$, then $U \cap E_1$ is a neighborhood of $x$ in $E_1$. So the claim follows from the already proven statement.
%\end{proof}
%
%\begin{remark}
%The continuity of the paths of $\mathbf M$, proved in Lemma \ref{LemmaOmega1P1} holds in the sense of sequentially continuity, see \cite{Doh05}. Since $([0,\infty),|\cdot|)$ is a metric space, the mapping $\mathbf X_s : [0,\infty) \to E$ is continuous, $\mathbf P_x$-a.s. Moreover, $\mathbf P_x$-a.s. $X_s \in E_1$, so $X_s \in C([0,\infty),E_1)$.
%So if $u \in C(E_1)$ we have that the mapping $u \circ X_t : [0,\infty) \to \R$ is $\mathbf P_x$-a.s. continuous.
%\end{remark}

\section{The underlying elliptic gradient Dirichlet form} \label{Dirichlet_form}

This section is devoted to prove Proposition \ref{PropDirichletform} of the introduction. We assume Condition \ref{AssEllipConstruction} and fix a $p>d$ as in Condition \ref{AssEllipConstruction}. We recall that the underlying Dirichlet form, the closure of \eqref{EqDirichletFormIntroduction}, serves as standard example in \cite{FOT94} and \cite{MR92} and can also be constructed assuming much weaker conditions. The stronger conditions on the density $\varrho$ and the matrix $A$ are required to apply our general construction scheme which yields the existence of an associated $\mathcal L^p$-strong Feller process that solves the martingale problem pointwisely. The full strength of the conditions becomes clear in this section and in Section \ref{SecConseqRegularity}. As already mentioned in the introduction, in case $A=I$ we are back in the setting of \cite{AKR03}.

%\begin{Cond} \label{Ass minimal}
%Let $\Omega \subseteq \mathbb{R}^n$, $n \in \mathbb{N}$, be an open set. Let $\mu$ be a measure on $(\Omega, \mathcal{B}(\Omega))$ having density $\varrho: \Omega \rightarrow \mathbb{R}$ w.r.t.~the Lebesgue measure $dx$ on $\Omega$ that satisfies
%\begin{enumerate}
%\item [(i)]
%$\varrho \geq 0$ $dx$-a.e.~on $\Omega$, $\varrho \in L^{1}_{\text{loc}}(\Omega,dx)$ and $\int_{U}\varrho ~dx > 0$ for all open $U \subseteq \Omega$,
%\item[(ii)]
%$\varrho$ satisfies the Hamza condition, i.e., $\varrho=0~dx$-a.e.~on $\Omega\setminus \mathcal{R}_\varrho(\Omega)$, where
%\begin{align*}
%\mathcal{R}_\varrho(\Omega):=\left\{ x\in\Omega \, \Bigg| \, \int_{\{y\in\Omega \, |\, |x-y|\le\varepsilon\}}\varrho^{-1}(y)dy<\infty\mbox{ for some }\varepsilon>0\right\}.
%\end{align*} .
%\end{enumerate}
%We assume $A=(a_{ij}(\cdot))_{i,j=1}^n: \Omega \rightarrow \mathbb{R}^{n \times n}$ to be measurable, symmetric such that
%\begin{enumerate}
%\item [(iii)]
%$A$ is locally strictly elliptic $dx$-a.e.~on $\Omega$.
%\item [(iv)]
%$a_{ij} \in L^1_{loc}(\Omega)$ for all $1 \leq i,j \leq n$.
%\end{enumerate}
%\end{Cond}

%We define the symmetric, non-negative definite bilinear form $(\E,D)$ on $L^2(\Omega,\mu)$ by
%
%\begin{Df}
%\begin{align*}
%&\E(f,g) := \int_{\Omega} ~ \big(A \nabla f,\nabla g\big)_{2} ~ d\mu = \int_{\Omega} ~ \sum_{i,j=1}^n (\partial_i f\,a_{ij} \, \partial_j g) ~ d\mu,~ f,g \in D:=C_0^\infty(\Omega).
%\end{align*}
%\end{Df}
%
%Then we obtain the following
\begin{remark} \label{RmHolderContinuousVersionRho}
Due to Condition \ref{AssEllipConstruction}(i),(ii) we conclude $\varrho \in H^{1,p}_{\text{loc}}(\Omega,dx)$. This follows analogously as in \cite[Cor.~2.2]{AKR03} by replacing $\mathbb{R}^d$ through $\Omega$ in the latter. In particular, we can choose a unique H"older continuous $dx$-version of $\varrho$ in $\Omega$. We fix this version and denote it again by $\varrho$.
\end{remark}

The following properties of our Dirichlet form are well-known. We just list the references.

\begin{proposition} \label{Pp Dirichlet form}
The bilinear form $(\mathcal{E},\mathcal{D})$ is closable on $L^2(\Omega,\mu)$. Its closure is denoted by $\Dir$ and is a symmetric, strongly local, regular Dirichlet form on $L^2(\Omega,\mu)$.
\end{proposition}

\begin{proof}
Since $\varrho$ is continuous, $\varrho$ fulfills the Hamza condition, see \cite[Ch.~II]{MR92}. Hence closability follows as in \cite[Ch.~II,~Exe.~2.4]{MR92}. The Dirichlet property follows as in \cite[Ch.~II,~Exe.~2.7]{MR92}. Concerning the regularity statement, see \cite[Ch.~IV,~Sec.~4a)]{MR92}. It can easily be proven using the extended Stone Weierstrass theorem, see e.g.~\cite[Ch.~7,~Sec.~38]{Sim63}. The latter also shows that $(\E,\D)$ is indeed densely defined. Locality is proven in \cite[Ch.~V,~Exa.~1.12(i)]{MR92}. In an analogous way one proves the strong local property. This can easily be done using \cite[Theo.~3.1.2]{FOT94}, \cite[Prob.~3.1.1]{FOT94} in combination with the extended Stone Weierstrass theorem.
\end{proof}

We have the following proposition.

\begin{proposition}\label{Pp Generator}
Let $p$ be as in Condition \ref{AssEllipConstruction}. For $f \in C_c^\infty(\Omega)$ it holds $f \in D(L_p) \cap D(L_2)$ and
\begin{align*} 
L_p f = L_2 f = L f :=\sum_{i,j=1}^d\partial_i(a_{ij}\,\partial_j f) + \partial_i\left(\ln \varrho \right) a_{ij} \, \partial_j f 
\end{align*}
Rewriting $L f$, we have
\begin{align} \label{Generator}
Lf = \sum_{i,j=1}^d a_{ij}\,\partial_i \partial_j f + \sum_{i=1}^d b_i \, \partial_i f,
\end{align}
where $b_i := \sum_{j=1}^d \partial_j a_{ij} + \partial_j \left( \ln\varrho \right) a_{ij}$, $1 \leq i \leq d$.
\end{proposition}

\begin{proof}
Let $f \in C_c^\infty(\Omega)$. First we check that $L f$ is indeed an element of $L^p(\Omega,\mu) \cap L^2(\Omega,\mu)$: Clearly $\partial_i(a_{ij}\,\partial_j f) \in L^p(\Omega,\mu)$ since $f \in C_c^\infty(\Omega)$, $\partial_i a_{ij} \in H^{1,\infty}_{\text{loc}}(\Omega)$ and $\varrho \in L^1(\supp[f])$. Moreover, $a_{ij} \pdj f \in L^\infty(\supp[f])$. Together with the assumption $\nabla \ln (\varrho) \in L^p(\Om,\mu)$ this implies $\pdi(\ln \varrho) a_{ij} \pdj f \in L^p(\Om,\mu)$. Altogether, $Lf \in L^p(\Omega,\mu)$. Hence also $Lf \in L^2(\Om,\mu)$ because $Lf$ vanishes outside the compact set $\supp[f]$.\\
Now let $g \in C_c^\infty(\Omega)$. By assumption we have $\varrho \in H^{1,1}_{\text{loc}}(\Omega)$ and $a_{ij} \in H^{1,\infty}_{\text{loc}}(\Omega)$. Since $\supp[f]$ is compactly contained in $\Om$, this implies $u_{ij}:=\varrho\, a_{ij}\,\partial_j f\in H^{1,1}(\Omega)$. Using integration by parts we obtain
\begin{align*}
\E(f,g) &= \sum_{i,j=1}^d \int_\Omega \partial_i g  \left( a_{ij} \, \partial_j f \right)\,d\mu = \sum_{i,j=1}^d \int_\Omega \partial_i g \left( \varrho \, a_{ij} \, \partial_j f \right)\,dx = -\sum_{i,j=1}^d \int_\Omega g \, \partial_i(\varrho \,a_{ij} \, \partial_j f)\,dx\\
&= -\sum_{i,j=1}^d \int_\Omega g \left(\, \partial_i(a_{ij} \, \partial_j f) +  a_{ij} \, \frac{\partial_i (\varrho)}{\varrho}  \, \partial_j f \right) d\mu =(-Lf,g)_{L^2(\Omega,\mu)}
\end{align*}
and therefore also $\E(f,g) = (-Lf,g)_{L^2(\Omega,\mu)}$ for all $g \in D(\E)$. Using \cite[Prop.~2.16]{MR92}, we conclude $f \in D(L_2)$ and $L_2 f=Lf$. And since $L f \in L^p(\Om,\mu)$ we also have $f \in D(L_p)$ by definition of the $L^p$-generator.
%Hence $(1-L_2)|_{C^\infty_c(\Om)} = (1-\tilde{L})$. Since $f$ and $\tilde{L} f \in L^2(\Om,\mu) \cap L^p(\Om,\mu)$ there exists a unique $w_0 \in D(L_2) \cap D(L_p)$ such that $(1-L_p)w_0 = (1-L_2)w_0 = (1-\tilde{L})f \,= (1-L_2)f$.
%\mfootnote{Indeed. Set $\tilde f := (1-\tilde L)f \in L^2 \cap L^p$. Set $f_n := max(min(-n,\tilde f)n) 1_{|\tilde f| \ge \frac{1}{n}} \in L^1 \cap L^\infty$. By Lebesgue $f_n$ converges both in $L^p$ and $L^2$-norm to $\tilde f$. By construction of the $L^p$-resolvent as unique extension of $R^2_1|_{L^1 \cap L^\infty}$ it holds $R^2_1 f_n = R^p_1 f_n$. The first term converges to $R^2_1 \tilde f$, the lather converges to $R^p_1 \tilde f$, hence $R^2_1 \tilde f = R^p_1 \tilde f =: w_0$.}

%Hence $w_0 = f$, thus $f \in D(L_p)$ and $L_p w_0 = \tilde{L} f$.\\
\end{proof}

Next we prove that the set $\{ \varrho = 0 \}$ has capacity zero. This is essential to obtain that our associated $\mathcal L^p$-strong Feller process stays in $\{ \varrho > 0\}$ unless it reaches the cemetery. The following proof is a slight modification of the ones from \cite{Fuk85} and \cite{FG08}. Here the proof of the second mentioned paper also relies on the ideas from \cite{Fuk85}. 

\begin{proposition} \label{Pp Capacity}
$\{\varrho=0\}$ is of capacity zero w.r.t.~$\Dir$.
\end{proposition}

\begin{proof}
For each $i \in \mathbb{N}$, choose open and bounded subsets $G_i \subseteq \Omega$ satisfying $\overline{G_i} \subseteq G_{i+1}$ and $\Omega = \bigcup_{i \in \mathbb{N}}{G_i}$. So $\{ \varrho=0 \} = \bigcup_{i \in \mathbb{N}}{G_i \cap \{ \varrho=0 \} }$ and thus $\text{cap}_{\mathcal E}(\{ \varrho=0 \}) = \sup_{i \in \mathbb{N}} \, {\text{cap}_{\mathcal E}(G_i \cap \{ \varrho=0 \})}$. So we only need to show $\text{cap}_{\mathcal E}(G_i \cap \{ \varrho=0 \})=0$ for all $i \in \mathbb{N}$. Therefore, let $i \in \mathbb{N}$ and choose $f \in C_c^\infty(\Omega)$ with the property that $f$ is strictly positive on $G_i$ and vanishes outside $G_{i+1}$. For each $\varepsilon > 0$ define
\begin{align*}
\psi := \sqrt{\varrho},~ \psi_{\varepsilon}:= \psi \vee \varepsilon,~f_\varepsilon := (\ln{\psi_{\varepsilon}}) \cdot f.
\end{align*}
By Condition \ref{AssEllipConstruction}(i) we have $\psi  \in H^{1,2}_{\text{loc}}(\Omega)$, thus $\psi_{\varepsilon} \in H^{1,2}_{\text{loc}}(\Omega)$ and $f_\varepsilon \in H^{1,2}(\Omega)$. Next we claim  $f_\varepsilon \in \Dom$ and 
\begin{align} \label{form_E_1(varepsilon)}
\E_1(f_\varepsilon,f_\varepsilon) \leq \int_{G_{i+1}} {(A \nabla f_\varepsilon, \nabla f_\varepsilon)_{\text{euc}}\, d\mu} + \int_{G_{i+1}} {f_\varepsilon^2 \,d\mu}.
\end{align}
Indeed, choose $g_k \in C^\infty(\Omega)$, $k \in \mathbb{N}$, such that $g_k \rightarrow f_\varepsilon$ in $H^{1,2}(\Omega,dx)$ as $k \rightarrow \infty$. We may assume that all $g_k$ and $f_\varepsilon$ have support inside some compact set $K$, $K \subseteq G_{i+1}$, by multiplying all $g_k$ and $f_\varepsilon$ with some cut off function $\chi \in C_c^\infty(G_{i+1})$ satisfying $\chi=1$ on $\text{supp}[f]$. Then also $g_k {\rightarrow} f_\varepsilon$ and $\nabla g_k {\rightarrow} \nabla f_\varepsilon$ in $L^2(\Omega,\mu)$ as $k \rightarrow \infty$ since $\varrho$ is bounded on $K$ by continuity. Furthermore, let $M < \infty$ be the $L^\infty$-bound of the matrix $A$ on $K$. Thus $\E(g_k,g_k) \leq M \int_K {|\nabla g_k|^2 ~d\mu}$, $k \in \mathbb{N}$, and hence $\sup_{k \in \mathbb{N}} \, \E(g_k,g_k) < \infty$. Then \cite[Ch.~1,~Lem.~2.12]{MR92} implies  $f_\varepsilon \in \Dom$ and
\begin{align*}
\E(f_\varepsilon, f_\varepsilon) \leq \liminf_{k \rightarrow \infty} \E(g_k,g_k) = \int_{G_{i+1}} {(A \nabla f_\varepsilon, \nabla f_\varepsilon)_{2}~d\mu},
\end{align*}
where the last equality holds because $A \,\nabla g_k {\rightarrow} A \,\nabla f_\varepsilon$ in $L^2(\Omega,\mu)$ as $k \rightarrow \infty$. Hence \eqref{form_E_1(varepsilon)} follows. In particular,
\begin{align} \label{eq_capacity_estimate}
\E_1(f_\varepsilon,f_\varepsilon) \leq M_{i} \int_{G_{i+1}} {( \nabla f_\varepsilon, \nabla f_\varepsilon)_{2}\, d\mu} + \int_{G_{i+1}} {f_\varepsilon^2 \,d\mu}.
\end{align}
with $M_{i} < \infty$ being the $L^\infty$-bound of $A$ on $G_{i+1}$. Now the right hand side of \eqref{eq_capacity_estimate} can be computed analogously as in \cite{Fuk85} and converges to some finite constant $C_{i}$ as $\varepsilon \rightarrow 0$. Here $C_{i}$ only depends on $G_{i+1}$. Let $\lambda > 0$. The calculation from \cite{Fuk85} shows 
\begin{align*}
\text{cap}_{\mathcal E} (G_i \cap \{ \varrho=0 \}) \leq \limsup_{n \in \mathbb{N}} \bigl\{ \text{cap}_{\mathcal E}(\{ |f_{\frac{1}{n}}(\cdot)| > \lambda \}) \bigr\} \leq \frac{1}{\lambda^2} \limsup_{n \in \mathbb{N}} \E_1(f_{\frac{1}{n}},f_{\frac{1}{n}}) = \frac{1}{\lambda^2} C_{i}.
\end{align*}
So $\lambda \rightarrow \infty$ implies the claim.
\end{proof}

\section{An elliptic regularity result} \label{regularity_results}

In this section we prove an elliptic regularity result, see Theorem \ref{Thm Elliptic regularity result} below. This provides the desired regularity Condition \ref{CondRegularityDomain}(i) for the Dirichlet form defined in \eqref{EqDirichletFormIntroduction} and is applied mainly in Theorem \ref{Thm Strong Feller Resolvent}.
\medskip

A similar version of the regularity result below is remarked in \cite[Rem.~2.15]{BKR01}, see also \cite[Cor.~2.13]{BKR01}. The assumptions on the coefficient functions therein are too strong for applying it in our case. Thus we need a generalized version of the statement from \cite[Rem.~2.15]{BKR01}. To obtain our desired regularity result we make use of the ideas from \cite[Theo.~2.8]{BKR01}. Our result then also is a generalization of \cite[Cor.~2.10]{BKR01}. We note that most parts of the proof concerning the following theorem are identical to the proof of \cite[Theo.~2.8]{BKR01}. Thus we only indicate the differences.
\medskip

In the sequel we adopt the notation from \cite{BKR01}: $A=(a_{ij})_{i,j=1}^d$ is always assumed to be a symmetric, nonnegative, $d\times d$-matrix valued measurable mapping on an open subset $U \subseteq \mathbb{R}^d$. For $\kappa \geq 0$ we write $A \geq \kappa I$ $dx$-a.e.~if $\kappa \sum_{i=1}^d\xi_i^2\ \leq \sum_{i,j=1}^d a_{ij}(x)\xi_i\xi_j$ holds for all $\xi \in \mathbb{R}^d$ and for $dx$-a.e.~$x \in {U}$. $A$ is called locally strictly elliptic $dx$-a.e.~on $U$ if for each $K \subseteq U$, $K$ compact, there exists some $\kappa_K > 0$ such that $A \geq \kappa_K I$ holds $dx$-a.e.~on $K$.

%\begin{Pp} \label{Pp1}
%Let $U \subseteq \mathbb{R}^n$ be open. Let $\nu$ and $\gamma$ be two locally finite (signed) Borel measures on $U$ and let $a_{ij}, b_i, c$ be locally $\nu$-integrable. Assume that
%\begin{align*}
%\int_{U} ( \sum_{i,j=1}^n a_{ij} \, \partial_i \partial_j \varphi + \sum_{i=1}^n b_i \, \partial_i \varphi + c \, \varphi ) ~d\nu = \int_{U} \varphi ~d\gamma \text{ for all nonnegative } \varphi \in C_0^\infty(U).
%\end{align*}
%If $A=(a_{ij})$ is locally H"older continuous and non-degenerate (i.e.,~$A(x)$ is non-singular for all $x \in U$), then $\nu$ has a density, which belongs to $L^r_{\text{loc}}(U,dx)$ for every $r \in [1,n'[$, where $\frac{1}{n'} + \frac{1}{n} =1$.
%\end{Pp}
%
%\begin{proof} See \cite [Corollary 2.3]{BKR01}.
%\end{proof}
%
%\begin{Pp} \label{Pp2}
%Let $p>n\geq2$, $q \in ]1,\infty[$, $R_1>0$. Let $a_{ij} \in H^{1,p}_{\text{loc}}(B_{R_1})$ be continuous and assume $A \geq \kappa I$ for some
%$\kappa > 0$. Then there exists $R_0 \in [0,R_1]$ with the following property: If $R<R_0$ and $\gamma$ is a (signed) measure of finite total variation on $B_R$ such that for any $\varphi \in C_0^2(B_R)$ we have
%\begin{align*}
%\left|\int_{B_R} \, \sum_{i,j=1}^n a_{ij} \, (\partial_i \partial_j \varphi) \, d\gamma \right| \leq N \left|\left| \nabla \varphi \right|\right|_{L^q(B_R)},
%\end{align*}
%where $N \in ]0,\infty[$ is independent of $\varphi$, then $\gamma \in H^{1, q' \wedge p}_{0}(B_R)$ (where $\gamma$ is identified with its density).
%\end{Pp}
%
%\begin{proof} See \cite [Theorem 2.7]{BKR01}.
%\end{proof}
\begin{theorem} \label{Thm Elliptic regularity result}
Let $U \subseteq \mathbb{R}^d$ be open with $d \geq 2$. Let $\nu$ be a locally finite (signed) Borel measure on $U$ that is absolutely continuous w.r.t.~the Lebesgue measure $dx$ on $U$. Let $a_{ij} \in H^{1,p}_{\text{loc}}(U)$ for some $p>d$ and assume $A$ to be locally strictly elliptic $dx$-a.e.~on $U$. Let either $b_i, c \in L^p_{\text{loc}}(U,dx)$ or $b_i, c \in L^p_{\text{loc}}(U,\nu)$ and let $f \in L^p_{\text{loc}}(U,dx)$. Assume that one has 
\begin{align}\label{Elliptic1}
\int_{U} \big( \sum_{i,j=1}^d a_{ij} \, \partial_i \partial_j \varphi + \sum_{i=1}^d b_i \, \partial_i \varphi + c \, \varphi \big) ~d\nu = \int_{U} f \, \varphi~dx ~~ \mbox{for all} ~~ \varphi \in C_c^\infty(U),
%\mfootnote{If $b_i, c$ are locally $dx$-integrable, then there are also locally $\nu$-integrable since $\nu$ is locally finite.} 
\end{align} 
where we assume that $b_i, c$ are locally $\nu$-integrable. Then $\nu$ has a density in $H^{1,p}_{\text{loc}}(U)$ that is locally H"older continuous.
\end{theorem}

\begin{proof} First let $\beta:= \left| b \right| + \left|c\right| +1$. We either have $\beta \in L^p_{\text{loc}}(U,dx)$ or $\beta \in L^p_{\text{loc}}(U,\nu)$. \\
\textit{Step 1}. By Sobolev embedding $A$ has a locally H"older continuous version. From now on, we fix this version and denote it with the same symbol. Continuity now implies that $A$ is locally strictly elliptic everywhere on $U$. In particular, $A$ is non-degenerate. Then \cite[Cor.~2.3]{BKR01} implies that $\nu$ has a density in $L^r_{\text{loc}}(U,dx)$ for every $r \in (1,d')$, where $d' := \frac{d}{d-1}$ is the dual exponent.\newline
\textit{Step 2}. We have $p>d$, i.e., $p'<d'$. So choose $r\in (p',d')$. Define $q=q(r):=\frac{pr}{pr-p-r} > 1$ and $q'=q'(r)=\frac{pr}{p+r}>1$. Analogously as in \cite[Theo.~2.8]{BKR01} one shows $\beta \nu \in L^{q'}_{\text{loc}}(U,dx)$. Further note that $f \in L^{q'}_{\text{loc}}(U,dx)$ since $q' \leq p$. Now choose an arbitrary $x_0 \in U$ and $R>0$ with $\overline{B_R} \subset U$ where ${B_R}=B_R(x_0)$. Let $\eta \in C_c^\infty(B_R)$. By using \eqref{Elliptic1}, a similar calculation as in \cite[Theo.~2.8]{BKR01} yields
\begin{align} \label{Elliptic2}
\left|\int_{B_R} \, \sum_{i,j=1}^d a_{ij} \, (\partial_i \partial_j \varphi) \eta \nu ~dx\right| \leq N \left|\left| \nabla \varphi \right|\right|_{L^q(B_R)},~\varphi \in C^2_0(B_R),
\end{align}
where $N < \infty$ is a constant not depending on $\varphi$. Here $C^2_0(B_R) := C^2(\overline{B_R}) \cap \{ u \, | \, u_{|{\partial B_R}}=0 \}$.\\
To be more precise, we follow the lines of $(2.22)$ in the proof of \cite[Theo.~2.8]{BKR01} but perform the estimate
\begin{align*}
\left|\int_{B_R} \, \sum_{i,j=1}^n a_{ij} \, \partial_i \partial_j (\varphi \eta) \,\nu ~dx\right| &\leq N_1 \int_{B_R} \big( \left( \left| \nabla \varphi \right| +\left| \varphi \right| \right) \left( \left|\beta \nu \right| + \left| f \right| \right) \big)~dx \\
&\leq N_1 \left|\left| \nabla \varphi \right|\right|_{L^q(B_R)} \left\| \left|\beta \nu \right| + \left| f \right| \right\|_{L^{q'}(B_R)} =  N_2 \left|\left| \nabla \varphi \right|\right|_{L^q(B_R)},
\end{align*}
where $N_1 < \infty$ and $N_2 < \infty$ are constants independent of $\varphi$. Herein the first inequality follows by \eqref{Elliptic1} and the second one by H"older inequality in combination with Poincar\'e inequality. Thus by $(2.22)$ in the proof of \cite[Theo.~2.8]{BKR01}, \eqref{Elliptic2} follows. \newline
\textit{Step 3}. Based on inequality \eqref{Elliptic2}, the conclusion $\nu \in H^{1,p}_{\text{loc}}(U)$ can be proven with the same method as presented in \cite[Theo.~2.8]{BKR01}.
\end{proof}

\begin{remark}
The statement of Theorem \ref{Thm Elliptic regularity result}(i) in case $a_{ij} \in C^\infty(U)$ is also contained in \cite{BKR97}, see \cite[Lem.~3,~Rem.~4(iii)]{BKR97}. For details see also \cite{Hen08}.
\end{remark}

Additionally, we need Morrey`s a priori estimate for which the original proof can be found in \cite[Thm.~5.5.5']{MR66}. For another detailed workout, see \cite{Sh06}.

\begin{theorem} \label{Thm Morrey}
Assume $p > d \geq 2$. Let $U\subset\mathbb{R}^d$ be open and bounded with $C^1$-boundary and let $V\subset\mathbb{R}^d$ be open and bounded such that $\overline{U}\subseteq V$. Let ${b}: V \to \mathbb{R}^d$ and ${c,e}: V \to \mathbb{R}$ such that
\begin{align*}
b_i \in L^p(V,dx) \mbox{ and } c,e \in L^q(V,dx)\mbox{ for }q:=\frac{dp}{d+p}>1.
\end{align*}
Let $a_{ij}:V \to \mathbb{R}$ be continuous for all $1 \leq i,j \leq d$. Furthermore, let $A \geq \kappa I$ for some $\kappa > 0$. Assume that $u \in H^{1,p}(U)$ is a solution of
\begin{align*}
\int_{U} \sum_{i=1}^d \big( \partial_i \varphi \, (\, \sum_{j=1}^d {a_{ij} \partial_j u} + b_i u )\big)+ \varphi (cu+e) ~ dx = 0\,\mbox{ for all } \, \varphi \in C_0^1(U),
\end{align*}
where $C_0^1(U) := C^1(\overline{U}) \cap \{ u \, | \, u_{\partial U} = 0 \}$.
Then we obtain for some constant $C < \infty$, independent of $e$ and $u$, the estimate
\begin{align*} 
\Vert u \Vert_{H^{1,p}(U)} \leq C \, ( \Vert e \Vert_{L^q(U,dx)} + \Vert u \Vert_{L^1(U,dx)}).
\end{align*}

\end{theorem}

\section{Construction of $\mathcal L^p$-strong Feller elliptic diffusions} \label{SecConseqRegularity}

As in \cite{AKR03} we apply the regularity result from the previous section to obtain estimates for the resolvent of $\Dir$, see Theorem \ref{Thm Strong Feller Resolvent}. The latter then directly allows us to use our general construction scheme and thus we can finally prove Theorem \ref{Th Resolvent Continuous} and \ref{DiffProcess} from the introduction. In this way, we can construct a $\mathcal L^p$-strong Feller process associated to our Dirichlet form $\Dir$. In this section we assume the same conditions as in Section \ref{Dirichlet_form} with the same $p>d$. Recall that $\varrho$ is chosen to be continuous, see Remark \ref{RmHolderContinuousVersionRho}.

\begin{remark} 
We denote by $(G^r_\lambda)_{\lambda > 0}$ the  strongly continuous sub-Markovian resolvent on $L^r(\Omega,\mu)$ associated to $(T_t^r)_{t >0}$, $r \in[1,\infty)$. Since $G^r_\lambda$ is the Laplace transform of $(T^r_t)_{t>0}$ the restriction of ${G^r_\lambda}$, $\lambda >0$, on $L^1(\Omega,\mu) \cap L^{\infty}(\Omega,\mu)$ coincides with $G^2_\lambda$. 
\end{remark}

Before proving the desired theorem we need one more Lemma.

\begin{lemma} \label{Lm Strong Feller}
Let $V$ be open and bounded with $\overline{V} \subset \{ \varrho > 0 \}$. Then we have $\hat{b}_i \in L^p(V,dx)$, where $\hat{b}_i := \sum_{j=1}^d{\partial_j a_{ij}} - b_i$, $1 \leq i \leq d$, and $b_i$ is defined as in Proposition \ref{Pp Generator}.
\end{lemma}

\begin{proof}
We have $\hat{b}_i \in L^p_{\text{loc}}(\Omega,\mu)$ and hence $\hat{b}_i \in L^p(V,\mu)$, i.e.,~ $\varrho |\hat{b}_i|^p \in L^1(V,dx)$. But since $\overline{V} \subseteq \{ \varrho >0 \}$ is compact and $\varrho$ is continuous, we conclude that $\varrho$ is bounded from below on $\overline{V}$ by some $\varrho_0 >0$. So $ |\hat{b}_i|^p \leq \frac{\varrho}{\varrho_0} |\hat{b}_i|^p$ on $V$ and therefore $|\hat{b}_i|^p \in L^1(V,dx)$.
\end{proof}

The following theorem is analogous to \cite[Cor.~2.3]{AKR03} for $A=Id$. For $A$ smooth see also \cite{Hen08}.
\begin{theorem} \label{Thm Strong Feller Resolvent}
Let $p$ be as in Condition \ref{AssEllipConstruction} and $\lambda>0$. Let $f\in L^p(\Omega,\mu)$. Then
\begin{align} \label{Equation Strong Feller Resolvent 1}
\varrho G^{p}_{\lambda} f \in H^{1,p}_{\text{loc}}(\{ \varrho > 0 \})
\end{align}
and for any open ball $B \subseteq \overline{B} \subseteq \{ \varrho > 0 \}$ there exists a finite constant $C_5$  independent of $f$ such that
\begin{align}\label{Equation Strong Feller Resolvent 2}
\Vert \varrho G^{p}_\lambda f \Vert_{H^{1,p}(B)} \leq C_5 \big( \Vert G^{p}_\lambda
f \Vert_{L^1(B,\mu)}+ \Vert f \Vert_{L^p(B,\mu)} \big).
\end{align}
\end{theorem}

\begin{proof}
Let us first assume $f \in C^{\infty}_{c}(\Omega) \subseteq L^1(\Omega,\mu) \cap L^{\infty}(\Omega,\mu)$. Then  $G^{p}_{\lambda} f = G^2_{\lambda} f \in L^{\infty}(\Omega,\mu) \cap L^2(\Omega,\mu)$. Let $\varphi \in C^{\infty}_{c}(\Omega)$. Symmetry of $L_2$ on $L^2(\Omega,\mu)$ leads to
\begin{align*}
{\left(\left(\lambda - L_2\right)\varphi,G^{p}_{\lambda} f\right)}_{L^2(\Omega,\mu)} &= 
{\left(\varphi,\left(\lambda - L_2\right) G^{2}_{\lambda} f\right)}_{L^2(\Omega,\mu)} 
= {\left(\varphi,f\right)}_{L^2(\Omega,\mu)}.
\end{align*}
Note that $G^{p}_{\lambda} f  \in L^1(\Omega,\mu)$. We define the locally finite signed Borel measure $\nu$ by \newline $\nu := G^{p}_{\lambda} f \, \mu$. Then $\nu$ is absolutely continuous w.r.t~$dx$. Using the representation of $L_2 \varphi$ (see \eqref{Generator}) we conclude
\begin{align} \label{***}
\int_{\Omega}{\Big( \sum_{i,j=1}^d a_{ij}\,\partial_i \partial_j \varphi + \sum_{i=1}^d b_i \, \partial_i \varphi -\lambda \varphi \Big)~ d\nu} = \int_{\Omega}{ \left( - f \varrho \right) \varphi ~dx}.
\end{align}
We plan to apply Theorem \ref{Thm Elliptic regularity result} and therefore have to check all the necessary assumptions: The conditions on $A=(a_{ij})$ are clearly satisfied. So $G^{p}_{\lambda} f \in L^{\infty}(\Omega,\mu)$ together with $b_i \in L^p_{\text{loc}}(\Omega,\mu)$ implies $b_i \in L^p_{\text{loc}}(\Omega,\nu)$. Finally, clearly $ \lambda \in L^p_{\text{loc}}(\Omega,\nu),~ f \varrho \in C(\Omega) \subseteq L^p_{\text{loc}}(\Omega,dx).$ \\
Thus \ref{Thm Elliptic regularity result} implies $\varrho G^{p}_{\lambda} f \in H^{1,p}_{\text{loc}}(\Omega,dx)$, in particular, (\ref{Equation Strong Feller Resolvent 1}) is shown for $f \in C^{\infty}_{c}(\Omega)$. Define $u:= \varrho  G^{p}_{\lambda} f$.  Using integration by parts,  $u$ solves
\begin{align*}
\int_{\Omega}{\sum_{i=1}^d {\Big( \partial_i \varphi \, ( \sum_{j=1}^d {a_{ij} \partial_j u} + \hat{b}_i u )\Big)} + \varphi (\lambda u - \varrho f) ~dx} = 0,
\end{align*}
where $\hat{b}_i = \sum_{j=1}^d{\partial_j a_{ij}} - b_i$. Now let $B$ be any open ball in $\mathbb{R}^d$ with $B \subseteq \overline{B} \subseteq \{ \varrho>0 \}$. Then by Lemma \ref{Lm Strong Feller} and Theorem  \ref{Thm Morrey} we get
\begin {align*}
\Vert \varrho G^{p}_{\lambda} f \Vert_{H^{1,p}(B)} \leq  C_6 \, \left( \Vert f \varrho \Vert_{L^q(B,dx)} + \Vert \varrho G^{p}_{\lambda} f \Vert_{L^1(B,dx)} \right)
\end{align*}  
for some $C_6 < \infty$ independent of $f$. Since $\Vert \cdot \Vert_{L^q(B,dx)}$ can be estimated from above by $\Vert \cdot \Vert_{L^p(B,dx)}$ and $|\varrho|^{p-1}$ is bounded on $B$, there exists $C_5 < \infty$ independent of $f$ such that
\begin {align} \label{****}
\Vert \varrho G^{p}_{\lambda} f \Vert_{H^{1,p}(B)} \leq  C_5 \, \left( \Vert f \Vert_{L^p(B,\mu)} + \Vert G^{p}_{\lambda} f \Vert_{L^1(B,\mu)} \right).
\end{align} 
So (\ref{Equation Strong Feller Resolvent 2}) is shown for $f \in C^{\infty}_{c}(\Omega)$. Now let $f \in L^p(\Omega,\mu)$ arbitrary. Choose $f_n \in C^{\infty}_{c}(\Omega)$ with $f_n \rightarrow f$ in $L^p(\Omega,\mu)$ as $n\rightarrow \infty$. Due to $\mu(B) < \infty$ we obtain $G^{p}_{\lambda} f_n \rightarrow G^{p}_{\lambda} f$ in $L^1(B,\mu)$ as $n\rightarrow \infty$. Hence by (\ref{****}) the sequence $\varrho G^{p}_{\lambda} f_n$  converges to some $g \in H^{1,p}(B)$ as $n \rightarrow \infty$. It is easy to see that $g$ coincides with $\varrho G^{r}_{\lambda} f$ $dx$-a.e.~on $B$. Thus $\varrho G^{r}_{\lambda} f \in H^{1,p}(B)$ and (\ref{Equation Strong Feller Resolvent 2}) is shown for $f \in L^p(\Omega,\mu)$. 
\end{proof}

Now we can prove Theorem \ref{Th Resolvent Continuous} from the introduction. We use that if $f,g \in C^{0,\beta}(\overline{U})$, $U \subset \mathbb{R}^d$ being open and bounded, then $fg \in C^{0,\beta}(\overline{U})$ and it holds 
\begin{align} \label{eq_proof_theorem_resolvent_continuous}
\|fg\|_{C^{0,\beta}(\overline{U})} \leq \|f\|_{C^{0,\beta}(\overline{U})} \|g\|_{C^{0,\beta}(\overline{U})}.
\end{align}

\begin{proof}[\textbf{Proof of Theorem \ref{Th Resolvent Continuous}}]
Let $x \in \{ \varrho > 0\}$, $r>0$ such that $\overline{B_r(x)} \subset \{ \varrho > 0\}$. We write $B = B_r(x)$. By Theorem \ref{Thm Strong Feller Resolvent} we have $\varrho G^p_\lambda f \in H^{1,p}(B)$. Using Sobolev's embedding we obtain that $\varrho G^p_\lambda f$ has a  H"older continuous version of index $\beta$ in $\overline{B}$. We denote the (unique) continuous version of $\varrho G^p_\lambda f$ on $\{ \varrho > 0\}$ by $\widetilde{\varrho G^p_\lambda f}$. Also note that $\varrho \in C^{0,\beta}(\overline{B})$, see Remark \ref{RmHolderContinuousVersionRho}.  Furthermore, $\varrho$ is bounded from below on $\overline{B}$ by some $\varrho_0 > 0$. So by using the global Lipschitz continuity of the map $\varphi(x):= x^{-1}$, $x \in [\varrho_0, \infty)$, we conclude $\varrho^{-1} \in {C^{0,\beta}(\overline{B})}$. Hence $G^p_\lambda f$ admits a  H"older continuous version of index $\beta$ in $\overline{B}$. Thus there exists a continuous $\mu$-version of $G^p_\lambda f$ on $\{\varrho > 0\}$. This version is unique, since $\mu$ is strictly positive on open sets. Together with $D(L_p) = \mathcal R (G^p_\lambda)$ we get the embedding $D(L_p) \emb C(\{ \varrho > 0 \})$.
Finally, inequality \eqref{Estimate Th Resolvent Continuous} follows by using \eqref{eq_proof_theorem_resolvent_continuous}, the estimate from Sobolev's embedding and estimate $\eqref{Equation Strong Feller Resolvent 2}$ from Theorem \ref{Thm Strong Feller Resolvent}.
\end{proof}

Finally, we have everything at hand to construct the $\mathcal L^p$-strong Feller process associated to $\Dir$, i.e., we prove Theorem \ref{DiffProcess} from the introduction.

\begin{proof}[\textbf{Proof of Theorem \ref{DiffProcess}}]
We plan to apply Theorem \ref{TheoProcess}. Set $E:=\Om$, $d$ the restriction of the Euclidean metric to $\Om$, $E_1 :=  \{ \varrho > 0\}$. Clearly $(E,d)$ is a locally compact separable metric space. Furthermore, set $\mu := \varrho dx$. Then $\mu$ is finite on compact sets and the assumption $\varrho > 0$ $dx$-a.e.~implies that $\mu$ has full topological support. By Proposition \ref{PropDirichletform} we have that $(\mathcal E,D(\mathcal E))$ is a strongly local, regular, symmetric Dirichlet form. It is left to check that the additional regularity conditions from \ref{CondRegularityDomain} are satisfied. Indeed, we have $\text{cap}_{\mathcal E}(\{ \varrho = 0 \})=0$ by Proposition \ref{PropDirichletform} and $C^\infty_c(\Om) \subset D(L_p)$. Clearly $C^\infty_c(\Om)$ is point separating in $\{ \varrho > 0 \}$ in the required sense. From Theorem \ref{Th Resolvent Continuous} we get that the embedding $D(L_p) \emb C(E_1)$ exists and is locally continuous. So all assumptions of Theorem \ref{TheoProcess} are fulfilled and hence we obtain a diffusion process $\mathbf M$ satisfying the stated properties except of the strong Feller property for the resolvent. The latter property can be shown with the help of Theorem \ref{ClassicalStrongFeller}. Indeed, let $(u_n)_{n \in \N}$ be a sequence in $D(L_p)$ such that $(f_n)_{n \in \N}$, $f_n:=(1-L_p)u_n \in L^p(\Om,\mu)$, $n \in \mathbb{N}$, is a bounded sequence in $L^\infty(\Om,\mu)$. Since $u_n = R^p_1 f_n$ and $R^p_1$ is sub-Markovian, we conclude that $(u_n)_{n \in \N}$ is also bounded in $L^\infty(\Om,\mu)$. Now for each $x \in \{ \varrho > 0\}$ we find an $r>0$ such that $\overline{B_r(x)} \subset \{ \varrho > 0\}$. Note that estimate \eqref{Estimate Th Resolvent Continuous} from Theorem \ref{Th Resolvent Continuous} implies
\begin{align*}
\Vert \widetilde{u_n} \Vert_{C^{0,\beta}(\overline{B_r(x)})} \leq K_1  (\Vert u_n \Vert_{L^\infty(\Omega,\mu)} + \Vert f_n \Vert_{L^\infty(\Omega,\mu)}),
\end{align*}
where $K_1$ is a finite constant independent of all $u_n$ and $f_n$, $n \in \mathbb{N}$. From this we get that $(\widetilde{u_n})_{n \in \N}$ is equicontinuous in $x$. Finally, Theorem \ref{ClassicalStrongFeller} implies that the associated resolvent of kernels $(R_\lambda)_{\lambda > 0}$ are even strong Feller. This finishes the proof.
\end{proof}

\textbf{\underline{Acknowledgement:}} We thank Torben Fattler and Sven-Ole Henning for preliminary work. Financial support by the DFG through the project GR 1809/8-1 is gratefully acknowledged.

%\bibliographystyle{alpha}
%\bibliography{../../../Bibliography/processes,../../../Bibliography/semigroups,../../../Bibliography/elliptic_pde,../../../Bibliography/workouts,../../../Bibliography/functional_analysis,../../../Bibliography/probability,../../../Bibliography/topology}

\hfill
\begin{minipage}[b]{10cm}
B. Baur\\
Department of Mathematics\\
University of Kaiserslautern\\
P.O. Box 3049, 67653 Kaiserslautern, Germany.\\
E-mail: baur@mathematik.uni-kl.de\\
URL: \url{http://www.mathematik.uni-kl.de/~wwwfktn}

\bigskip
M. Grothaus\\
Department of Mathematics\\
University of Kaiserslautern\\
P.O. Box 3049, 67653 Kaiserslautern, Germany.\\
E-mail: grothaus@mathematik.uni-kl.de\\
URL: \url{http://www.mathematik.uni-kl.de/~grothaus}

\bigskip
P. Stilgenbauer\\
Department of Mathematics\\
University of Kaiserslautern\\
P.O. Box 3049, 67653 Kaiserslautern, Germany.\\
E-mail: stilgenb@mathematik.uni-kl.de\\
URL: \url{http://www.mathematik.uni-kl.de/~wwwfktn}

\end{minipage}

\end{document}